\documentclass
{imsart}

\RequirePackage{amsthm,amsmath,amsfonts,amssymb}
\RequirePackage[numbers,sort&compress]{natbib}
\RequirePackage[colorlinks,citecolor=blue,urlcolor=blue]{hyperref}
\RequirePackage{graphicx}

\startlocaldefs
\theoremstyle{definition}

\newtheorem{theorem}{Theorem}

\newtheorem{proposition}{Proposition}
\newtheorem{corollary}{Corollary}

\theoremstyle{definition}
\newtheorem{definition}{Definition}

\numberwithin{equation}{section}

\endlocaldefs

\usepackage{natbib}
\bibpunct[:]{(}{)}{,}{a}{}{,}

\usepackage{subcaption}

\theoremstyle{definition}

\newcommand{\argmin}{\mathop{\rm arg~min}\limits}

\begin{document}

\begin{frontmatter}
\title{Local Fr\'echet Regression via RKHS embedding and Its Applications to Data Analysis on Manifolds}
\runtitle{LFR via RKHS embedding}

\begin{aug}
\author[A]{\fnms{Iida}~\snm{Yuki}
\ead[label=e1]{cdghp0812@keio.jp}},
\author[A]{\fnms{Shiraishi}~\snm{Hiroshi}\ead[label=e2]{shiraishi@math.keio.ac.jp}}
\and \author[B]{\fnms{Ogata}~\snm{Hiroaki}\ead[label=e3]{hiroakiogata@tmu.ac.jp}}
\address[A]{Department of Mathematics, 
Keio University, Yokohama, Kanagawa, Japan, \\ \printead{e1,e2}}
\address[B]{Faculty of Economics and Business Administration,
Tokyo Metropolitan University, Hachioji, Tokyo, Japan\printead[presep={,\ }]{e3}}
\end{aug}

\begin{abstract}
Local Fr\'echet Regression (LFR) is a nonparametric regression method for settings in which the explanatory variable lies in a Euclidean space and the response variable lies in a metric space. It is used to estimate smooth trajectories in general metric spaces from noisy observations of random objects taking values in such spaces. Since metric spaces form a broad class of spaces that often lack algebraic structures such as addition or scalar multiplication—characteristics typical of vector spaces—the asymptotic theory for conventional random variables cannot be directly applied. As a result, deriving the asymptotic distribution of the LFR estimator is challenging.
In this paper, we first extend nonparametric regression models for real-valued responses to Hilbert spaces and derive the asymptotic distribution of the LFR estimator in a Hilbert space setting. Furthermore, we propose a new estimator based on the LFR estimator in a reproducing kernel Hilbert space (RKHS), by mapping data from a general metric space into an RKHS. Finally, we consider applications of the proposed method to data lying on manifolds and construct confidence regions in metric spaces based on the derived asymptotic distribution.
\end{abstract}


\end{frontmatter}


\section{Introduction}

The analysis of data with complex structures in non-Euclidean spaces has been gaining popularity across various domains of data science (see [\cite{MR3258095}] for an overview). A fundamental concept for describing the distribution of random elements in non-Euclidean spaces is the Fr\'echet mean [\cite{MR27464}], which generalizes the notion of the expectation of a random variable to a metric space. [\cite{MR3909947}] extended this concept to conditional distributions and developed least squares and nonparametric regression methods for settings where the explanatory variable lies in Euclidean space and the response variable lies in a metric space. These methods are referred to as Global Fr\'echet Regression (GFR) and Local Fr\'echet Regression (LFR), respectively.

LFR extends local linear regression [\cite{MR1193323}, \cite{MR1212173}], a nonparametric regression method for real-valued responses, to settings where the response variable lies in a metric space. This allows for the estimation of a smooth trajectory on the metric space from noisy observations of random objects taking values in that space. The LFR estimator is defined as the solution to a minimization problem that seeks a weighted Fr\'echet mean, where the weights are determined by a kernel function applied to each sample point. In general, this minimization problem is not analytically solvable, and even numerical solutions require the computation of distances between points in the metric space, resulting in high computational cost. Moreover, unlike in vector spaces, the complex structure of points in metric spaces prevents the direct application of conventional asymptotic theory for real-valued random variables. As a result, deriving the asymptotic distribution of the LFR estimator is challenging. Indeed, while [\cite{MR3909947}] and [\cite{MR4441132}] established convergence rates for the distance between the true regression function and the LFR estimator, they did not address its asymptotic distribution. Even in the case where the response space is a separable Hilbert space, only the asymptotic distribution of the GFR estimator has been derived, with no corresponding results for LFR.

In this paper, we extend the classical nonparametric regression model for real-valued responses to the Hilbert space setting by taking advantage of the fact that the LFR estimator admits an explicit form when the response space is a Hilbert space, and derive the asymptotic distribution of the LFR estimator in Hilbert spaces based on the asymptotic theory of local polynomial regression [\cite{MR1383587}]. Next, for random objects lying in metric spaces that are not Hilbert spaces, we propose mapping the data into a reproducing kernel Hilbert space (RKHS) using a positive definite kernel, following the kernel method framework [\cite{scholkopf2002learning}], and then performing LFR in the RKHS. The resulting LFR estimator in the RKHS is then mapped back to elements in the original metric space using the distance function, giving rise to a new estimator of the true regression function based on LFR. This approach allows for the estimation of regression trajectories in metric spaces even when the structure of the space is too complex to compute LFR directly, or when the computational cost is prohibitive.
[\cite{MR3909947}] conducted simulations of LFR using data distributed on the unit 2-sphere $S^2$ as an example. The unit sphere is one of the simplest Riemannian manifolds, and data distributed on it are known as directional data. Directional data are observations that represent directions or angles, and due to their periodic nature, conventional statistical methods for real valued data often fail to apply appropriately
(see [\cite{MR1828667}, \cite{Jammalamadaka2001Topics}] for examples). 
Nonparametric estimation methods for directional data include kernel density estimation [\cite{MR919843}, \cite{MR2772221}], smoothing splines on nonlinear manifolds [\cite{su2012fitting}], and kernel-based trend estimation for time series trends [\cite{MR3565479}].
Regression models where the response variable lies on a more general finite-dimensional manifold, including directional data, have also been widely studied. These include intrinsic models that directly utilize the geometric structure of the manifold [\cite{shi2009intrinsic}, \cite{niethammer2011geodesic}, \cite{MR3104017}, \cite{MR3611755}], and extrinsic models that embed the manifold into a Euclidean space [\cite{lin2017extrinsic}, \cite{lee2021robust}]. Moreover, extensions of the classical Nadaraya-Watson smoother to manifold-valued responses have been proposed [\cite{MR2214065}, \cite{davis2010population}, \cite{hinkle2012polynomial}, \cite{MR2965956}, \cite{kilin2023stability}].
In this paper, we analyze random objects on more complex manifolds - specifically, on the ellipsoid and the torus - rather than on the unit circle or unit sphere, which are commonly studied in directional statistics. With respect to data on ellipsoids, statistical analyses have been conducted in the context of astronomy and fluid dynamics, targeting various datasets defined on spheroidal surfaces, including celestial bodies, as reported in [\cite{MR3357394}, \cite{MR3720779}, \cite{wang2017ellipsoidal}, \cite{kilin2023stability}]. In addition, several studies, including [\cite{MR3834306},  \cite{MR4451377}, \cite{nodehi2025torus}, \cite{price2025data}], have formulated RNA structures, celestial orbits, vortex in the nonlinear Schr\"odinger equation, and so on as data residing on the two-dimensional torus, and have developed a range of statistical methodologies for their analysis. In this study, we demonstrate the utility of LFR and our proposed method for data on the ellipsoid and the torus through simulation experiments, and we also construct an approximate confidence region based on the asymptotic distribution.

The structure of this paper is as follows. In Section \ref{Sec2}, we introduce the Local Fr\'echet Regression model and the LFR estimator as defined in
[\cite{MR3909947}]. In Section \ref{Sec3}, we consider a nonparametric regression model in a Hilbert space and derive the asymptotic distribution of the LFR estimator. In Section \ref{Sec4}, we propose a new estimator by mapping data in a general (non-Hilbert) metric space to an RKHS, performing LFR there, and then mapping the RKHS-based LFR estimator back to the original metric space. We also discuss the construction of approximate confidence regions based on the asymptotic distribution. In Section \ref{Sec5}, we present simulation studies to illustrate applications of the proposed methods to data on manifolds. Section \ref{Sec6} concludes. 
Proofs of the theoretical results presented in Section \ref{Sec3} and \ref{Sec4} are provided in the \hyperref[app]{Appendix}.
Throughout this paper, we use the symbol $\leadsto$ to denote weak convergence.

\section{Local Fr\'{e}chet Regression}\label{Sec2}
Let $(\Omega, d)$ be a metric space. Consider a random pair $(X, Y) \sim F$, where $X$ takes values in $\mathbb{R}$ and $Y$ takes values in $\Omega$, and $F$ is the joint distribution of $(X, Y)$ on $\mathbb{R} \times \Omega$. Denote the marginal distributions of $X$ and $Y$ by $F_X$ and $F_Y$, respectively, and assume that the mean $\mu = \mathrm{E}[X]$ and variance $\Sigma = \mathrm{Var}(X) > 0$ exist. We also assume that the conditional distributions $F_{X|Y}$ and $F_{Y|X}$ exist.
In this general setting, $Y$ is referred to as a random object. In the classical regression problem where $\Omega = \mathbb{R}$, the goal is to estimate the conditional mean $m(x) = \mathrm{E}[Y|X = x]$.
This regression function can be characterized as $m(x) = \arg\min_{y \in \mathbb{R}} \mathrm{E}[d_E^2(Y, y)|X = x]$,
where $d_E$ denotes the Euclidean distance.
[\cite{MR3909947}] introduced the Fr\'echet regression function by generalizing this representation: replacing the Euclidean distance $d_E$ with a general metric $d$, and replacing $\mathbb{R}$ with a general metric space $\Omega$.
\begin{definition}[{\textit{Fr\'{e}chet regression function} [\cite{MR3909947}]}]\label{de1}
\begin{align}
m_{\oplus}(x) = \argmin_{\omega \in \Omega} M_{\oplus} (\omega, x), \quad \quad
M_{\oplus}(\cdot, x) = \mathrm{E}[d^2(Y, \cdot) | X = x].  \label{2.1}
\end{align}
\end{definition}
A nonparametric regression model used to estimate $m_{\oplus}(x)$ is the \textit{local F\'echet regression model} 
proposed by [\cite{MR3909947}].
In the classical regression setting where $(\Omega, d) = (\mathbb{R}, d_E)$, [\cite{MR1383587}] consider the following nonparametric regression model:
\begin{align}
Y = m(X) + \sigma(X) \varepsilon, \label{2.2}
\end{align}
where $\sigma^2(x) = \mathrm{Var}(Y \mid X = x)$, and the random variable $\varepsilon$ is independent of $X$, with \\
$\mathrm{E}[\varepsilon] = 0$ and $\mathrm{Var}(\varepsilon) = 1$.
Given an i.i.d. sample $(X_1, Y_1), \ldots, (X_n, Y_n)$ from model \eqref{2.2}, the local linear (LL) estimator [\cite{MR1193323}, \cite{MR1212173}] of $m(x)$
is given by $\hat{m}_{\mathrm{LL}}(x) = \hat{\beta}_0$, where
\begin{align*}
(\hat{\beta}_0, \hat{\beta}_1) = \arg\min_{\beta_0, \beta_1 \in \mathbb{R}} \frac{1}{n} 
\sum_{i = 1}^{n} K_h(X_i - x) (Y_i - \beta_0 - \beta_1 (X_i - x))^2, 
\end{align*}
with a kernel function $K$, bandwidth $h$, and $K_h(\cdot) = h^{-1} K(\cdot / h)$.\\
{[\cite{MR3909947}]} interpreted $\hat{\beta}_0$ as a weighted Fr\'echet mean:
\begin{align}
\tilde{l}(x) &= \mathrm{E}[s(X, x, h) Y] = \arg\min_{y \in \mathbb{R}} \mathrm{E}[s(X, x, h) d_E^2 (Y, y)], \label{2.3}
\end{align}
where
$s(X, x, h) = \frac{1}{\sigma_0^2} K_h(X - x) \left\{\mu_2 - \mu_1 (X - x)\right\}$,\ $\mu_j = \mathrm{E}[K_h(X - x)(X - x)^j]$,\\
$(j = 0, 1, 2),\  \sigma_0^2 = \mu_0 \mu_2 - \mu_1^2$, 
and then extended the concept of local regression to random objects $Y \in \Omega$.
\begin{definition}[{\textit{Local Fr\'echet Regression model} [\cite{MR3909947}]}]\label{de2}
\begin{align}
\tilde{l}_{\oplus}(x) = \arg\min_{\omega \in \Omega} \tilde{L}_n(\omega), \quad 
\tilde{L}_n(\omega) = \mathrm{E}[s(X, x, h) d^2(Y, \omega)], \label{2.4}
\end{align}
where the dependency on $n$ is through the bandwidth sequence $h = h_n$.
\end{definition}

Given i.i.d. observations $(X_1, Y_1), \ldots, (X_n, Y_n) \overset{\mathrm{i.i.d.}}{\sim} F$, $\tilde{l}_{\oplus}(x)$ can be estimated using
\[
\hat{\mu}_j = \frac{1}{n} \sum_{i = 1}^{n} K_h(X_i - x) (X_i - x)^j, (j = 0, 1, 2), \quad \hat{\sigma}_0^2 = \hat{\mu}_0 \hat{\mu}_2 - \hat{\mu}_1^2,
\]
and the empirical weight function
\[
s_{in}(x, h) = \frac{1}{\hat{\sigma}_0^2} K_h(X_i - x) \left[\hat{\mu}_2 - \hat{\mu}_1 (X_i - x)\right].
\]

\begin{definition}[{\textit{LFR estimator} [\cite{MR3909947}]}]\label{de3}
\begin{align}
\hat{l}_{\oplus}(x) = \argmin_{\omega \in \Omega} \hat{L}_n(\omega), \quad 
\hat{L}_n(\omega) = \frac{1}{n} \sum_{i = 1}^{n} s_{in}(x, h) d^2(Y_i, \omega). \label{2.5}
\end{align}
\end{definition}
\section{Local Fr\'echet Regression on a Hilbert Space}\label{Sec3}
In general metric spaces, algebraic structures such as addition and scalar multiplication, which are defined in vector spaces, are not necessarily available. Moreover, because the Fr\'echet regression function $m_{\oplus}(\cdot)$ is defined without assuming a structural model like \eqref{2.2}, it is extremely difficult to derive the asymptotic distribution of the LFR estimator. Even in [\cite{MR3909947}], the discussion is limited to convergence rates depending on the metric. Furthermore, the minimization problems in \eqref{2.4} and \eqref{2.5} require the computation of distances between points in the metric space and generally cannot be solved analytically.
On the other hand, when $\Omega$ is a Hilbert space, the solution $\tilde{l}_{\oplus}(x)$ to the minimization problem exists uniquely, and the estimator $\hat{l}_{\oplus}(x)$ can be explicitly computed.
In the following, for Hilbert space $\mathcal{H}$, 
We denote by $\langle \cdot , \cdot \rangle_{\mathcal{H}}$ the inner product on $\mathcal{H}$. 
The norm is given by $\| \omega \|_{\mathcal{H}} := \langle \omega , \omega \rangle^{1/2}_{\mathcal{H}}$, and the distance between two points $x, y \in \mathcal{H}$ is defined as $d(x, y) := \| x - y \|_{\mathcal{H}} = \langle x - y, x - y \rangle^{1/2}_{\mathcal{H}}$. 
Assume further that $\mathbb{E}[\| Y \|_{\Omega}^2] < \infty$ when $\Omega$ is a Hilbert space (i.e., $\Omega = \mathcal{H}$).

\begin{proposition}\label{pr1}
Suppose that $\Omega$ is a Hilbert space and $\mathbb{E}[\| Y \|_{\Omega}^2] < \infty$. Then, for any
fixed $x \in \mathbb{R}$, the minimizer $\tilde{l}_{\oplus}(x)$ in \eqref{2.4} exists and is unique.
\end{proposition}

\begin{proposition}\label{pr2}
Suppose that $\Omega$ is a Hilbert space. Then, for any fixed $x \in \mathbb{R}$,
the minimizer $\hat{l}_{\oplus}(x)$ in \eqref{2.5} is explicitly given by
\begin{align}
\hat{l}_{\oplus}(x) = \frac{1}{n} \sum_{i=1}^{n} s_{in}(x, h) Y_i. \label{3.1}
\end{align}
\end{proposition}

We extend the nonparametric regression model \eqref{2.2} for real-valued response variables to the case where $\Omega$ is a general Hilbert space, and derive the asymptotic distribution of the LFR estimator.
\begin{definition}[\textit{Nonparametric Regression Model for Hilbert Space-Valued Response Variables}]\label{de4}
\begin{align}
Y = m_{\oplus}(X) + \sigma(X) \varepsilon, \label{3.2}
\end{align}
where $\varepsilon$ is a random element in Hilbert space $\Omega$, independent of $X$, satisfying $\mathrm{E}[\langle \varepsilon, \omega \rangle_{\Omega}] = 0$ for all $\omega \in \Omega$, and 
$\mathrm{E}[ \| \varepsilon \|_{\Omega}^2 ] = \mathrm{E}[ \langle \varepsilon, \varepsilon \rangle_{\Omega} ] = 1$. 
Also, $\sigma^2(x) = \mathrm{E}[ \| Y - m_{\oplus}(x) \|_{\Omega}^2 \mid X = x ]$. 
\end{definition}
Then, $m_{\oplus}(\cdot)$ in Definition \ref{de4} coincides with $m_{\oplus}(\cdot)$ in Definition \ref{de1} 
(see SUPPLEMENTARY MATERIAL \hyperref[S1]{S.1} for details).
The regression model \eqref{3.2} is an extension of the regression model \eqref{2.2} to Hilbert spaces. 
By considering model \eqref{3.2}, and based on the asymptotic theory of local polynomial regression~[\cite{MR1383587}], we can derive the asymptotic normality of $\hat{l}_{\oplus}(x)$ when $\Omega$ is a real Hilbert space. We impose the following assumptions for a fixed $x \in \mathbb{R}$.
\begin{description}
\item[(N1)]
$X$ has a density function $f$, and $f$ is continuous and positive in some neighborhood $N_x$ of $x$.

\item[(N2)]
The kernel $K$ is a symmetric probability density function with bounded support, and both $\int u^2 K(u) \, du$ and $\int K^2(u) \, du$ are finite.

\item[(N3)]
$h = h_n > 0$ satisfies $h \to 0$, $nh \to \infty$, and $nh^5 \to 0$ as $n \to \infty$.

\item[(O1)]
For all $\omega \in \Omega$, the function $m_\omega(\cdot) = \langle m_{\oplus}(\cdot), \omega \rangle_{\Omega}$ is twice continuously differentiable on $N_x$.

\item[(O2)]
For all $\omega \in \Omega$, $\sigma^2(\cdot) \, \mathrm{E}[\langle \varepsilon, \omega \rangle^2_{\Omega}]$ is continuous and positive on $N_x$.

\item[(O3)]
For all $\omega \in \Omega$, $\mathrm{E}[\langle Y, \omega \rangle^4_{\Omega} \mid X = \cdot]$ is bounded on $N_x$.
\end{description}
(N1)--(N3) are conditions on the density of $X$, the kernel, and the bandwidth, which are required for the asymptotic normality of the local linear (LL) estimator [\cite{MR1193323}, \cite{MR1212173}] for real-valued responses. 
(O1)--(O3) are the conditions on $m(\cdot)$ and $\sigma(\cdot)$ in model \eqref{2.2}, adapted to the model \eqref{3.2}.
\begin{theorem}\label{th1}
Let $\Omega$ be a real Hilbert space, and suppose that a random sample $(X_1, Y_1), \ldots, (X_n, Y_n)$ 
is drawn from the model \eqref{3.2}. 
If (N1)--(N3) and (O1)--(O3) hold, then for any $\omega \in \Omega$, the following holds:
\begin{align*}
& \sqrt{nh} \left\{  \left \langle  \hat{l}_\oplus(x) - m_\oplus(x) , \omega \right \rangle _\Omega -  
 \left( \frac{\int u^2 K(u) \, du}{2} \cdot \frac{\partial ^2}{\partial x^2} \left \langle m_\oplus(x) , \omega \right \rangle _\Omega \right) h^2 \right\} \\
 {} \leadsto & \; \mathcal{N} \left( 0, \; \frac{\int K^2(u) \, du}{f(x)} \sigma^2(x) \, \mathrm{E}[\langle \varepsilon, \omega \rangle_\Omega^2] \right).
\end{align*}

In particular, if we take $\omega = m_\oplus(x)$, we obtain the following expression :
\begin{align*}
& \sqrt{nh} \left\{ \left \langle  \hat{l}_\oplus(x), m_\oplus(x) \right \rangle_\Omega - \| m_\oplus(x) \|_\Omega^2 
- \left( \frac{\int u^2 K(u) \, du}{2} \cdot\left.  \frac{\partial ^2}{\partial x^2} 
\left \langle m_\oplus(x) , \omega \right  \rangle _\Omega \right|_{\omega = m_{\oplus}(x)} \right) h^2 \right\} \\
{} \leadsto & \; \mathcal{N} \left( 0, \; \frac{\int K^2(u) \, du}{f(x)} \sigma^2(x) \, \mathrm{E}[\langle \varepsilon, m_\oplus(x) \rangle_\Omega^2] \right).
\end{align*}\end{theorem}
If $\Omega$ is separable, then there exists a countable orthonormal basis $\{ e_i \}_{i = 1}^{\infty} \subset \Omega$. In this case, the difference 
$\hat{l}_\oplus(x) - m_\oplus(x)$ can be represented as:
\[
\hat{l}_\oplus(x) - m_\oplus(x) = \sum_{i = 1}^{\infty} a_{in}(x, X_1, \ldots, X_n) e_i, \quad a_{in}(x, X_1, \ldots, X_n) \in \mathbb{R}.
\]
Taking $\omega = e_j$, we obtain
\[
\left\langle \hat{l}_\oplus(x) - m_\oplus(x), e_j \right\rangle_{\Omega} = a_{jn}(x, X_1, \ldots, X_n),
\]
and hence the asymptotic normality of the $j$-th coefficient $a_{jn}(x, X_1, \ldots, X_n)$ of $\hat{l}_\oplus(x) - m_\oplus(x)$ is obtained.
\mbox{}\\

Next, using the result of Theorem \ref{th1}, we consider constructing a confidence region for $m_{\oplus}(x)$ on Hilbert space.
To derive the convergence rate of the estimator $\hat{l}_{\oplus}(x)$, we make the following assumptions similarly to [\cite{MR3909947}]. For simplicity, we assume that the marginal density $f$ of $X$, within the joint distribution $F$, has unbounded support, and 
consider a fixed point $x \in \mathbb{R}$ such that $f(x) > 0$.
\begin{description}
\item[(P1)] Let $B_\delta(m_{\oplus}(x)) \subset \Omega$ be the ball of radius $\delta$ centered at $m_{\oplus}(x)$ and
$N(\eta, B_\delta(m_{\oplus}(x)), d)$ be its covering number using balls of size $\eta$. Then
\begin{align*}
\int_{0}^{1} \sqrt{1 + \log N(\delta \eta, B_\delta(m_{\oplus}(x)), d)} \, d\eta = O(1) \ \  \quad \text{as } \delta \to 0.
\end{align*}

\item[(K0)] The kernel $K$ is a probability density function symmetric around zero. Furthermore, defining $K_{ij} = \int_{\mathbb{R}} K^i(u) u^j \,du$,  
$|K_{14}|$, $|K_{44}|$, and $|K_{26}|$ are finite.

\item[(L0)] 
The object $m_{\oplus}(x)$ exists and is unique. For all $n$, $\tilde{l}_{\oplus}(x)$ and $\hat{l}_{\oplus}(x)$ exist and are unique, 
the latter almost surely. Additionally, for any $\varepsilon > 0$, 
\[
\inf_{d(\omega, m_{\oplus}(x)) > \varepsilon} \left \{ M_{\oplus}(\omega, x) - M_{\oplus}(m_{\oplus}(x), x) \right \} > 0, 
\]
\[
\liminf_{n} \inf_{d(\omega, \tilde{l}_{\oplus}(x)) > \varepsilon} \left \{ \tilde{L}_n(\omega) - \tilde{L}_n(\tilde{l}_{\oplus}(x)) \right \} > 0.
\]

\item[(L1)] The marginal density $f$ of $X$, as well as the conditional densities $g_y$ of $X|Y = y$, exist and are twice continuously differentiable,
the latter for all $y \in \Omega$, and $\sup_{x,y} |g''_y(x)| < \infty$. Additionally, for any open $U \subset \Omega$, $\int_U dF_{Y|X}(y|x)$ is continuous 
as a function of $x$.

\item[(L2)]
There exists $\eta_1 > 0, C_1 > 0$ and $1 < \beta_1 < 2$ such that
\[
M_{\oplus}(\omega, x) - M_{\oplus}(m_{\oplus}(x), x) \geq C_1d(\omega, m_{\oplus}(x))^{\beta_1}, 
\]
provided $d(\omega, m_{\oplus}(x)) < \eta_1$.
\item[(L3)]
There exists $\eta_2 > 0, C_2 > 0$ and $1 < \beta_2 < 2$ such that
\[
\liminf_{n} \left\{ \tilde{L}_n(\omega) - \tilde{L}_n(\tilde{l}_{\oplus}(x))  \right\} \geq C_2d(\omega, \tilde{l}_{\oplus}(x))^{\beta_2},
\]
provided $d(\omega, m_{\oplus}(x)) < \eta_2$.

\end{description}
(P1) and (L0) are based on empirical process theory and are used to control the behavior of the objective function near the minimizer.
(K0) and (L1) are standard assumptions in local regression estimation.
(L2) and (L3) provide the convergence rates between $m_{\oplus}(x)$ and $\tilde{l}_{\oplus}(x)$, and between $\tilde{l}_{\oplus}(x)$ and $\hat{l}_{\oplus}(x)$, respectively. 
Since the convergence is with respect to $(nh)^{1/2}$, we additionally assume $\beta_1, \beta_2 < 2$.

Fix some $\alpha \in (0,1)$ and $m_{\oplus}(x) \in \Omega$, and define
\begin{align*}
    &B_x = \left( \frac{\int u^2 K(u) \, du}{2} \cdot
    \left. \frac{\partial^2}{\partial x^2} 
    \left\langle m_{\oplus}(x), \omega \right\rangle_{\Omega} 
    \right|_{\omega = m_{\oplus}(x)} \right) h^2 \quad \text{(asymptotic bias)}, \\
    &V_x = \frac{\int K^2(u) \, du}{f(x)} \sigma^2(x) \, \mathrm{E}[\left\langle \varepsilon, m_{\oplus}(x) \right\rangle_{\Omega}^2]
    \quad \text{(asymptotic variance)}.
\end{align*}

\begin{corollary}\label{cr1}
Assume the conditions of Theorem \ref{th1}, and suppose that consistent estimators $\hat{B}_x$ and $\hat{V}_x$ of $B_x$ and $V_x$, respectively, satisfy
\[
\hat{B}_x = B_x + o_P \left((nh)^{-1/2}\right), \quad \hat{V}_x = V_x + o_P (1).
\]
Additionally, assume conditions (P1), (K0), and (L0)--(L3) hold. Then, define
\begin{align*}
CI_{\alpha} &:= \left\{ \omega \in \Omega \;\middle|\; 
-2\hat{B}_x - 2\sqrt{\frac{\hat{V}_x}{nh}}\, z_{\alpha/2}
\le 
\left\| \omega - \hat{l}_{\oplus}(x) \right\|_{\Omega}^2
\le 
-2\hat{B}_x + 2\sqrt{\frac{\hat{V}_x}{nh}}\, z_{\alpha/2}
\right\},
\end{align*}
it follows that
\[
P\left( m_{\oplus}(x) \in CI_{\alpha} \right) \to 1 - \alpha \quad \text{as } n \to \infty.
\]
\end{corollary}
By applying the result of Corollary \ref{cr1}, we obtain a confidence region for $m_{\oplus}(x)$ in the form of a spherical region centered at $\hat{l}_{\oplus}(x)$, when $\Omega$ is a real Hilbert space.
\mbox{}\\

Nonparametric estimation of $m_\oplus(\cdot)$ for response variables taking values in a general metric space $(\Omega, d)$ has predominantly been conducted using the Nadaraya-Watson-type estimator 
[\cite{steinke2008non}, \cite{hein2009robust}, \cite{davis2010population}, \cite{MR2736019}], which is defined as:
\begin{align}
\hat{m}^{\mathrm{NW}}_\oplus(x) = \argmin_{\omega \in \Omega} \frac{1}{n} \sum_{i = 1}^{n} K_h(X_i - x) d^2(Y_i, \omega). \label{3.3}
\end{align}

If $\Omega$ is a Hilbert space, the minimizer in \eqref{3.3} can be obtained explicitly, just as with the LFR estimator. For an i.i.d. sample $(X_1, Y_1), \ldots, (X_n, Y_n) \overset{\mathrm{i.i.d.}}{\sim} F$, we have
\begin{align}
\hat{m}^{\mathrm{NW}}_\oplus(x) 
&= \argmin_{\omega \in \Omega} \frac{1}{n} \sum_{i = 1}^{n} K_h(X_i - x) \| Y_i - \omega \|_\Omega^2 
= \frac{\sum_{i = 1}^{n} K_h(X_i - x) Y_i}{\sum_{j = 1}^{n} K_h(X_j - x)}. \label{3.4}
\end{align}

Given an arbitrary $\omega_0 \in \Omega$, we consider two estimators of $\langle m_\oplus(x), \omega_0 \rangle_\Omega$, namely, 
$\langle \hat{l}_\oplus(x), \omega_0 \rangle_\Omega$ and 
$\langle \hat{m}^{\mathrm{NW}}_\oplus(x), \omega_0 \rangle_\Omega$. 
We compare their conditional biases and variances, as summarized in Table~\ref{table1}. 
(See SUPPLEMENTARY MATERIAL \hyperref[S2]{S.2} for detailed derivations.)

\begin{table}
\centering
\caption{Conditional bias and variance of estimators for $\langle m_\oplus(x), \omega_0 \rangle$: Comparison between local Fr\'echet regression (LFR) and Nadaraya-Watson (NW) estimators.}
\label{table1}
\begin{tabular}{lcc}
\hline
\textbf{Method} & \textbf{Bias} & \textbf{Variance} \\ \hline\hline
\textbf{NW}   
& $\frac{\int u^2 K(u) du}{2}  \left \{  \frac{\partial ^2}{\partial x^2} \left \langle m_\oplus(x) , \omega_0   \right  \rangle _\Omega + 
2 \frac{f'(x) \frac{\partial}{\partial x} \left \langle m_\oplus(x) , \omega_0   \right  \rangle _\Omega } {f(x)} \right \} h^2$
& $\frac{\int K^2(u) du}{f(x) nh} \sigma ^2 (x)  \mathrm{E}[\langle \varepsilon, \omega_0 \rangle _\Omega ^2]  $  
\\ & $ \qquad \qquad \qquad \qquad \qquad \qquad \qquad \qquad \qquad \qquad \qquad \quad+ o_{P}(h^2)$ & $\qquad \qquad \qquad \qquad + o_{P}((nh)^{-1})$
\\  \hline
\textbf{LFR}    
& $\frac{\int u^2 K(u) du}{2} \cdot \frac{\partial ^2}{\partial x^2} \left \langle m_\oplus(x) , \omega_0   \right  \rangle _\Omega h^2$
& $\frac{\int K^2(u) du}{f(x) nh} \sigma ^2 (x)  \mathrm{E}[\langle \varepsilon, \omega_0 \rangle _\Omega^2]  $  
\\ & $\qquad \qquad \qquad \qquad \qquad \qquad  + o_{P}(h^2)$ & $\qquad \qquad \qquad \qquad + o_{P}((nh)^{-1})$
  \\ \hline
\end{tabular}
\end{table}

As shown in Table~\ref{table1}, the bias terms of both NW and LFR estimators are of the same order, 
but the NW estimator includes an additional term:
$\int u^2 K(u) \, du \cdot \frac{f'(x)}{f(x)} \cdot \frac{\partial}{\partial x} \langle m_\oplus(x), \omega_0 \rangle_\Omega$,
which does not appear in the bias of the LFR estimator.
As for the variance, the leading terms are identical for both estimators.
In LFR, the bias depends on the second derivative $\frac{\partial^2}{\partial x^2} \langle m_\oplus(x), \omega_0 \rangle_\Omega$, whereas in NW, the bias also depends on the first derivative $\frac{\partial}{\partial x} \langle m_\oplus(x), \omega_0 \rangle_\Omega$, as well as $f'(x)$ and $f(x)$. 
Therefore, the absolute value of the bias for the NW estimator tends to be larger in regions where the density $f(x)$ is small 
or where the derivatives of $\langle m_\oplus(\cdot), \omega_0 \rangle_\Omega$ or $f$ are large.
\section{Local Fr\'echet Regression via RKHS embedding}\label{Sec4}
In Section \ref{Sec3}, we showed that when $\Omega$ is a Hilbert space, the LFR estimator can be obtained explicitly, and its asymptotic distribution can be derived using the inner product.
In this section, based on the concept of kernel methods[\cite{scholkopf2002learning}], we consider mapping data lying in a metric space $\Omega$ (which is not necessarily a Hilbert space) into a Reproducing Kernel Hilbert Space (RKHS) using a positive definite kernel, and then performing LFR in the RKHS. Furthermore, by mapping the obtained LFR estimator in the RKHS back to the original metric space using the distance function, we obtain a new estimator of the true regression function $m_{\oplus}(\cdot)$ based on LFR in the RKHS. This approach enables regression analysis even for data in metric spaces where distance computations are intractable and direct LFR estimation is not feasible.

Let $k : \Omega \times \Omega \to \mathbb{R}$ be a positive definite kernel on $\Omega$, and let $H_k$ denote the RKHS associated with the kernel $k$. Define the feature map $\Phi : \Omega \rightarrow H_k$ by $\Phi(\omega) = k(\cdot, \omega)$, i.e., $\Phi$ maps a point $\omega \in \Omega$ to the function $k(\cdot, \omega) \in H_k$. For any $f \in H_k$ and $\omega \in \Omega$, the identity 
$\langle f, k(\cdot, \omega) \rangle_{H_k} = f(\omega)$ holds, which is known as the reproducing property.
We assume the following condition:
\begin{align*}
\mathrm{E}[ \| \Phi(Y) \|_{H_k}^2 ] = \mathrm{E}[ k(Y, Y) ] < \infty.
\end{align*}
By applying the map $\Phi$ to the data $Y_i \in \Omega$, $i = 1, \ldots, n$, we obtain the transformed data in the RKHS as
\[
\{ \Phi(Y_i) \}_{i = 1}^{n} = \{ k(\cdot , Y_i) \}_{i = 1}^{n}.
\]
Given the transformed response variable $\Phi(Y) \in H_k$ and a real-valued explanatory variable $X \in \mathbb{R}$, we consider the following nonparametric regression model.
\begin{definition}[\textit{Nonparametric Regression Model for Responses in $H_k$}]\label{de5}
\begin{align}
\Phi(Y) = m_{\oplus}^{H_k}(X) + \sigma(X) \varepsilon, \label{4.1}
\end{align}
where $\varepsilon$ is a random element in $H_k$, independent of $X$, satisfying
$\mathrm{E}[\langle \varepsilon, f \rangle_{H_k}] = 0$ for all $f \in H_k$, and
$\mathrm{E}[\| \varepsilon \|^2_{H_k}] = \mathrm{E}[\langle \varepsilon, \varepsilon \rangle_{H_k}] = 1.$
Also, $\sigma^2(x) = \mathrm{E}[\| \Phi(Y) - m_{\oplus}^{H_k}(X) \|^2_{H_k} | X = x]$.
\end{definition}

The function $m_{\oplus}^{H_k}(X) \in H_k$ is the Fr\'echet regression function for explanatory variable $X \in \mathbb{R}$ and response variable $\Phi(Y) \in H_k$. If the positive definite kernel $k$ satisfies
\begin{align}
\argmin_{f \in H_k} \mathrm{E}[\| \Phi(Y) - f \|^2_{H_k} | X = x] = \Phi(m_{\oplus}(x)), \label{4.2}
\end{align}
then we have $m_{\oplus}^{H_k}(x) = \Phi(m_{\oplus}(x))$.
Equation \eqref{4.2} holds when the mapping $\Phi$ preserves the ordering of distances on $\Omega$, i.e., for  $x_1, y_1, x_2, y_2 \in \Omega$,
$d(x_1, y_1) < d(x_2, y_2) \Rightarrow \| \Phi(x_1) - \Phi(y_1) \|_{H_k} < \| \Phi(x_2) - \Phi(y_2) \|_{H_k}$,
and when the noise in $Y$ around the true regression function $m_{\oplus}(x)$ is small. In what follows, we consider such kernels $k$.
From model \eqref{4.1}, given a random sample $(X_1, \Phi(Y_1)), \ldots, (X_n, \Phi(Y_n))$, we consider the LFR estimator $\hat{l}_{\oplus}^{H_k}(x)$. Since $H_k$ is a Hilbert space, it follows from Proposition \ref{pr2} that
\begin{align}
\hat{l}_{\oplus}^{H_k}(x) = n^{-1} \sum_{i = 1}^{n} s_{in}(x, h) \Phi(Y_i). \label{4.3}
\end{align}
We associate $\hat{l}_{\oplus}^{H_k}(x) \in H_k$ with a point in $\Omega$ as the solution to the following minimization problem, and define it as $\hat{l}_{\oplus}^{k}(x)$:
\begin{align}
\hat{l}_{\oplus}^{k}(x) &= \argmin_{\omega \in \Omega} \| \hat{l}_{\oplus}^{H_k}(x) - \Phi(\omega) \|^2_{H_k} \label{4.4} \\ 
&= \argmin_{\omega \in \Omega} \Big\{ k(\omega, \omega) - 2n^{-1} \sum_{i = 1}^{n} s_{in}(x, h) k(Y_i, \omega) \notag \\
& \qquad \qquad \qquad \qquad \quad+ n^{-2} \sum_{i = 1}^{n} \sum_{j = 1}^{n} s_{in}(x, h) s_{jn}(x, h) k(Y_i, Y_j) \Big\}. \notag
\end{align}
The transformation from the first line to the second and third in \eqref{4.4} uses the reproducing property of the kernel, i.e., $\forall x, y \in \Omega$, $\langle \Phi(x), \Phi(y) \rangle_{H_k} = k(x, y)$.
By transforming data $Y_i \in \Omega$ into $H_k$ via $\Phi$, we can consider local Fr\'echet regression in $H_k$. The resulting LFR estimator $\hat{l}_{\oplus}^{H_k}(x)$ is then mapped back to $\Omega$ using \eqref{4.4}, yielding an estimator $\hat{l}_{\oplus}^{k}(x)$ for $m_{\oplus}(x)$.
We refer to the estimator \eqref{4.4}, which is based on the positive definite kernel $k$, as the LFR$^k$ estimator.
Since equation \eqref{4.4} does not involve the distance $d$ on the metric space $\Omega$, it enables obtaining an LFR-based estimator on $\Omega$ even in cases where the structure of $\Omega$ is complex and solving the minimization problem \eqref{2.5} directly is intractable.
For simplicity, we assume that the marginal density $f$ of $X$, within the joint distribution $F$, has unbounded support and consider a fixed point $x \in \mathbb{R}$ such that $f(x) > 0$.
To establish the consistency of $\hat{l}_{\oplus}^{k}(x)$ for $m_{\oplus}(x)$, we impose the following assumptions.
\begin{description}
\item[(S1)] The metric space $(\Omega, d)$ is compact.

\item[(S2)] There exists a constant $C > 0$ such that for all $\omega_1, \omega_2 \in \Omega$,
\[
\| \Phi(\omega_1) - \Phi(\omega_2) \|_{H_k} \leq C \, d(\omega_1, \omega_2).
\]

\item[(S3)] There exists $M > 0$ such that for all $\omega \in \Omega$,
\[
\| \Phi(\omega) \|_{H_k} = k(\omega, \omega)^{1/2} < M.
\]

\item[(S4)] The object $m_{\oplus}(x)$ exists and is unique. Additionally, defining\\
$M_{\oplus}^{H_k}(f, x) = \mathrm{E}[ \| \Phi(Y) - f \|_{H_k}^2 | X = x ]$,\ \ 
$\tilde{L}_n^{H_k}(f) = \mathrm{E}[ s(X, x, h) \| \Phi(Y) - f \|_{H_k}^2 ]$,\\
$\tilde{l}_{\oplus}^{H_k}(x) = \arg\min_{f \in H_k} \tilde{L}_n^{H_k}(f)$,
assume that for any $\varepsilon > 0$,
\[
\inf_{\| f - \Phi(m_{\oplus}(x)) \|_{H_k} > \varepsilon} \left\{ M_{\oplus}^{H_k}(f, x) - M_{\oplus}^{H_k}(\Phi(m_{\oplus}(x)), x) \right\} > 0,
\]
\[
\liminf_{n} \inf_{\| f - \tilde{l}_{\oplus}^{H_k}(x) \|_{H_k} > \varepsilon} \left\{ \tilde{L}_n^{H_k}(f) - \tilde{L}_n^{H_k}(\tilde{l}_{\oplus}^{H_k}(x)) \right\} > 0.
\]

\item[(S5)] Defining
$\tilde{L}_n^k(\omega) = \| \tilde{l}_{\oplus}^{H_k}(x) - \Phi(\omega) \|_{H_k}^2$, 
$\tilde{l}_{\oplus}^k(x) = \arg\min_{\omega \in \Omega} \tilde{L}_n^k(\omega)$,
assume that for all $n \in \mathbb{N}$, both $\tilde{l}_{\oplus}^k(x)$ and $\hat{l}_{\oplus}^k(x)$ exist
and are unique, the latter almost surely. Additionally, for any $\varepsilon > 0$,
\[
\liminf_n \inf_{d(\omega, \tilde{l}_{\oplus}^k(x)) > \varepsilon} \left\{ \tilde{L}_n^k(\omega) - \tilde{L}_n^k(\tilde{l}_{\oplus}^k(x)) \right\} > 0.
\]
\end{description}
Since any compact metric space is separable, assumption (S1) implies that $H_k$ is separable by Lemma 4.33 in [\cite{MR2450103}].  
Assumptions (S2) and (S3) ensure that the distance between two functions in $H_k$, obtained by mapping two elements from $\Omega$ via $\Phi$, can be controlled by the distance between those original points in $\Omega$, and that the norm of the mapped functions is uniformly bounded, respectively.  
Assumptions (S4) and (S5) are based on the theory of M-estimators, and they control the behavior of the objective functions $M_{\oplus}^{H_k}$, $\tilde{L}_n^{H_k}$, and $\tilde{L}_n^k$ around their respective minimizers. Under these assumptions, we obtain the following result:
\begin{theorem}\label{th2}
If (K0), (L1), and (S1)--(S5) hold, and if $h = h_n \to 0$ and $nh^4 \to \infty$. then
\begin{align*}
d(m_{\oplus}(x), \hat{l}_{\oplus}^{k}(x)) = o_{P}(1).
\end{align*}
\end{theorem}
Even when the metric space $\Omega$ is so complex that computing distances is difficult and LFR estimators cannot be obtained directly from \eqref{2.5}, one can always explicitly obtain the LFR estimator $\hat{l}_{\oplus}^{H_k}(x)$ in the Hilbert space $H_k$ by applying the mapping $\Phi$ to the data $Y_i \in \Omega$, $i = 1, \ldots, n$.  
Theorem \ref{th2} shows that, under appropriate assumptions, the LFR$^k$ estimator $\hat{l}_{\oplus}^{k}(x)$, which corresponds to $\hat{l}_{\oplus}^{H_k}(x)$ mapped back onto $\Omega$, serves as a consistent estimator for $m_{\oplus}(x)$.
\mbox{} \\

Next, we consider constructing a confidence region on $\Omega$ using the result of Theorem \ref{th1}.  
By the reproducing property of the positive definite kernel, for any $g \in H_k$, we have:
\begin{align*}
  \left \langle  \hat{l}_\oplus ^{H_k}(x) - m_\oplus ^{H_k}(x) ,  g  \right \rangle_{H_k} &=
\left \langle  n^{-1} \sum_{i = 1}^{n} s_{in}(x, h) \Phi(Y_i) -  \Phi(m_{\oplus}(x)), \ g \right \rangle _{H_k} \\
&= n^{-1} \sum_{i = 1}^{n} s_{in}(x, h) g(Y_i) -  g(m_{\oplus}(x)).
\end{align*}
We impose the following assumptions (Q1)--(Q3), which correspond to assumptions (O1)--(O3) in Theorem \ref{th1}.
\begin{description}
\item[(Q1)]
For all $g \in H_k$, the function $g(m_{\oplus}(\cdot))$ is twice continuously differentiable on $N_x$.
\item[(Q2)]
For all $g \in H_k$, the function $\sigma^2(\cdot) \mathrm{E}[\langle \varepsilon, g \rangle_{H_k}^2]$ is continuous and positive on $N_x$.
\item[(Q3)]
For all $g \in H_k$, the conditional fourth moment $\mathrm{E}[g(Y)^4 \mid X = \cdot]$ is bounded on $N_x$.
\end{description}
\begin{theorem}\label{th3}
Suppose that a random sample $(X_1, \Phi(Y_1)), \ldots, (X_n, \Phi(Y_n))$ is drawn from model \eqref{4.1}.  
If (N1)--(N3) and (Q1)--(Q3) hold, then for any $g \in H_k$, the following holds:
\begin{align*}
& \sqrt{nh} \Bigg\{ \left \langle \hat{l}_\oplus^{H_k}(x) - m_\oplus^{H_k}(x), g \right \rangle_{H_k}
- \left( \frac{\int u^2 K(u) \, du}{2} \cdot \frac{\partial^2}{\partial x^2} g(m_\oplus(x)) \right) h^2 \Bigg\} \\
= & \sqrt{nh} \Bigg\{ n^{-1} \sum_{i = 1}^{n} s_{in}(x, h) g(Y_i) - g(m_\oplus(x))
- \left( \frac{\int u^2 K(u) \, du}{2} \cdot \frac{\partial^2}{\partial x^2} g(m_\oplus(x)) \right) h^2 \Bigg\} \\
\leadsto & \ \mathcal{N} \left( 0, \ \frac{\int K^2(u) \, du}{f(x)}  \sigma^2(x) \, \mathrm{E}[\langle \varepsilon, g \rangle_{H_k}^2] \right).
\end{align*}
In particular, for $\omega \in \Omega$, taking $g = \Phi(\omega) = k(\cdot, \omega)$ gives:
\begin{align*}
& \sqrt{nh} \Bigg\{ n^{-1} \sum_{i = 1}^{n} s_{in}(x, h) k(\omega, Y_i) - k(\omega, m_{\oplus}(x))  - \left( \frac{\int u^2 K(u) \, du}{2} \cdot \frac{\partial^2}{\partial x^2} k(\omega, m_{\oplus}(x)) \right) h^2 \Bigg\} \\
\leadsto & \ \mathcal{N} \left( 0, \ \frac{\int K^2(u) \, du}{f(x)} \sigma^2(x) \, \mathrm{E}[\varepsilon(\omega)^2] \right).
\end{align*}
\end{theorem}

Fix some $\alpha \in (0,1)$ and $m_{\oplus}(x) \in \Omega$, and define
\begin{align*}
&B_x^{H_k} = \left( \frac{\int u^2 K(u) \, du}{2} \cdot
\left. \frac{\partial^2}{\partial x^2} k(\omega, m_{\oplus}(x)) \right|_{\omega = m_{\oplus}(x)} \right) h^2 
\quad \text{(asymptotic bias)}, \\
&V_x^{H_k} = \frac{\int K^2(u) \, du}{f(x)} \, \sigma^2(x) \, \mathrm{E}[\varepsilon(m_{\oplus}(x))^2]
\quad \text{(asymptotic variance)}.
\end{align*}

\begin{corollary}\label{cr2}
Assume the conditions of Theorems \ref{th2} and \ref{th3}, and suppose that consistent estimators 
$\hat{B}_x^{H_k}$ and $\hat{V}_x^{H_k}$ satisfy
\[
\hat{B}_x^{H_k} = B_x^{H_k} + o_P \left((nh)^{-1/2}\right), \quad 
\hat{V}_x^{H_k} = V_x^{H_k} + o_P (1).
\]
Additionally, assume the analogues of (P1), (K0), and (L0)--(L3) hold in the reproducing kernel Hilbert space induced by the feature map $\Phi: \Omega \to H_k$.
Then, define
\begin{align*}
CI_{\alpha} := \left\{ \omega \in \Omega \;\middle|\; 
-2\hat{B}_x^{H_k} - 2\sqrt{\frac{\hat{V}_x^{H_k}}{nh}} \, z_{\alpha/2}
\le 
\left\| \Phi(\omega) - \hat{l}_{\oplus}^{H_k}(x) \right\|_{H_k}^2
\le 
-2\hat{B}_x^{H_k} + 2\sqrt{\frac{\hat{V}_x^{H_k}}{nh}} \, z_{\alpha/2}
\right\},
\end{align*}
it follows that
\begin{align*}
P \left( m_{\oplus}(x) \in CI_{\alpha} \right) \to 1 - \alpha, \ \ 
P\left( \hat{l}_{\oplus}^{k}(x) \in CI_{\alpha} \right) \to 1 - \alpha \quad \text{as } n \to \infty.
\end{align*}
\end{corollary}

By applying the result of Corollary \ref{cr2}, we can construct a confidence region for $m_{\oplus}(x)$ in the original metric space $\Omega$ where the data are distributed.
\section{Application to Data on Manifolds}\label{Sec5}

In this section, we consider an application to data that lie on a manifold. For example, data such as wind directions or ocean current directions can be represented as points on the unit circle, while data such as earthquake epicenters can be represented as points on the unit sphere; these are examples of manifold-valued data. These types of data, where each observation represents a direction or angle, are referred to as directional data. Due to their inherent periodicity, standard statistic and methods, such as the sample mean for real-valued data, may not be meaningful (see [\cite{MR1828667}, \cite{Jammalamadaka2001Topics}] for details).
On the other hand, since LFR focuses on the distances between data points and estimates the Fr\'echet mean with kernel-based weighting, it is considered well-suited for the analysis of data on manifolds, including directional data.
In the following, we present two simulation studies of regression analysis with manifold-valued response variables, where the computation of geodesic distances is typically expensive.

\begin{figure}
\centering
\includegraphics[width = 5cm]{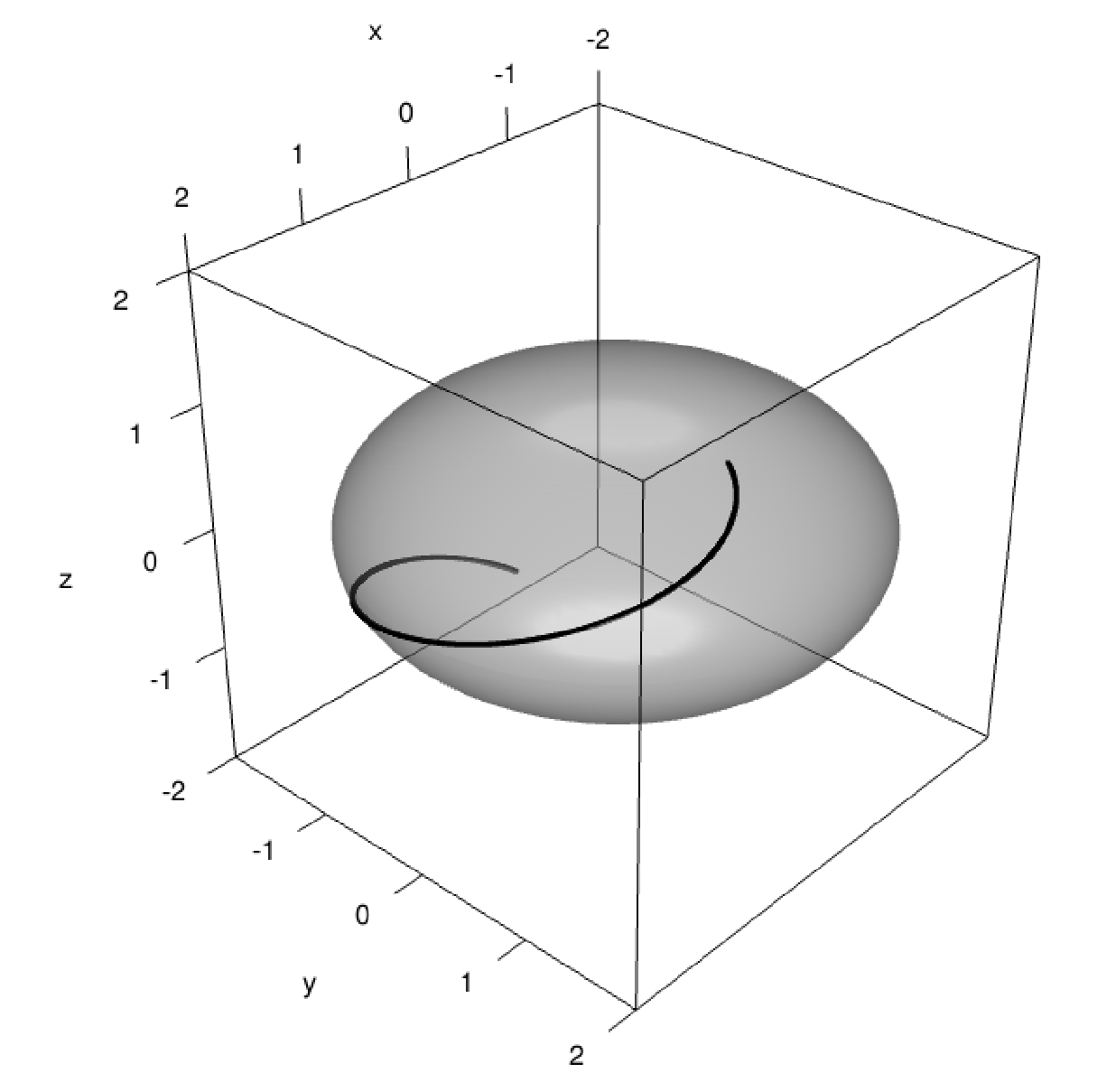} 
\caption{\textit{True regression curve \( m_{\oplus}(x) \) on \( E \), \( x \in (-0.5, 0.5) \).}}
\label{fig1}
\end{figure}

First, we consider the rotational ellipsoid $E$ as the metric space $\Omega$. 
The equation of a rotational ellipsoid with $z$ as the symmetry axis is given by
\begin{align}
\Omega = E := \left\{ (x, y, z)^{\top} \in \mathbb{R}^3 \mid \frac{x^2}{a^2} + \frac{y^2}{a^2} + \frac{z^2}{b^2} = 1 \right\}. \label{5.1}
\end{align}
The semi-axis $a$ is referred to as the \textit{equatorial radius} of the rotational ellipsoid, and the semi-axis $b$ is the distance from centre to pole along the symmetry axis. 
A point $(x, y, z)^{\top}$ on $E$ can be parameterized by the latitude $u \in \left[-\frac{\pi}{2}, \frac{\pi}{2}\right)$ and the longitude $v \in [-\pi, \pi)$ as follows:
\[
x = a \cos(u) \cos(v), \quad 
y = a \cos(u) \sin(v), \quad 
z = b \sin(u).
\]
We define the distance function $d$ on $E$ to be the geodesic distance. A geodesic on a surface (or more generally, on a Riemannian manifold) is a curve that gives the shortest path between two points, and the length of this curve is called the \textit{geodesic distance} (c.f., Chapter 6 in [\cite{MR1468735}]).
For two points $P_1, P_2 \in E$, the geodesic distance $d$ is defined using an auxiliary sphere. Let $\sigma$ denote the angle between $P_1$ and $P_2$ on the auxiliary sphere. Then the geodesic distance is given by:
\begin{align*}
d(P_1, P_2) := b \int_0^{\sigma} \sqrt{1 + k^2 \sin^2 \sigma'} \, d\sigma',
\end{align*}
where $k = e' \cos \alpha_0$, $e' = \sqrt{(a^2 - b^2)/b^2}$ is the second eccentricity of the $E$, and $\alpha_0$ is the azimuth of the geodesic at the equatorial intersection (c.f., \cite{rapp1993geometric}). 
In the special case where $a = b = 1$, i.e., when $E$ is the two-dimensional unit sphere $S^2$, we have $e' = 0$, and the geodesic distance reduces to an explicit formula: $d(P_1, P_2) = \sigma = \arccos(P_1^{\top} P_2)$.
On the other hand, on a rotational ellipsoid with $a \neq b$ (i.e., not a sphere), the geodesic equations cannot, in general, be solved analytically, 
therefore, computing the geodesic distance between two points requires numerical methods. Several algorithms have been developed for this purpose, such as the Vincenty method [\cite{vincenty1975direct}] and the Karney method [\cite{karney2013algorithms}], which achieve high accuracy but incur greater computational cost.
In the following, we perform a simulation study on the rotational ellipsoid $E$ with $a = 2$ and $b = 1$. 
We define the true regression function as a spiral curve on $E$, as shown in Fig.~\ref{fig1}:
\begin{align}
m_{\oplus}(x) = \left(2 \sqrt{1 - 3x^2} \cos(\pi x),\ 
                     2 \sqrt{1 - 3x^2} \sin(\pi x),\ 
                     \sqrt{3} x \right), 
\quad x \in (-0.5,\ 0.5). \label{5.2}
\end{align}

\begin{figure}
   \centering
  \begin{minipage}[b]{0.43\linewidth}
    \centering
    \includegraphics[keepaspectratio, scale=0.25]{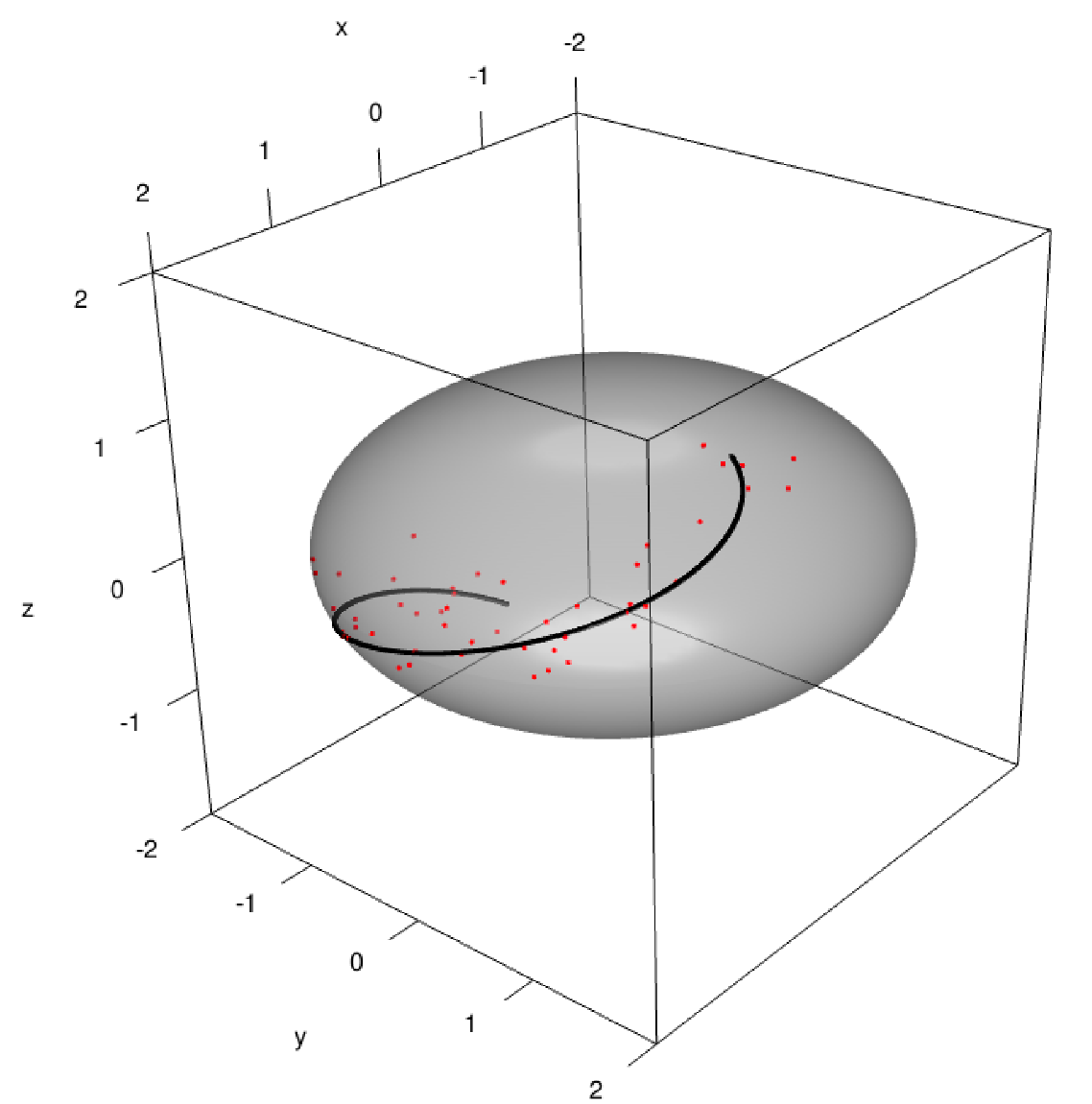}
    \subcaption{Low Noise, $n = 50$}
  \end{minipage}
  \begin{minipage}[b]{0.43\linewidth}
    \centering
    \includegraphics[keepaspectratio, scale=0.244]{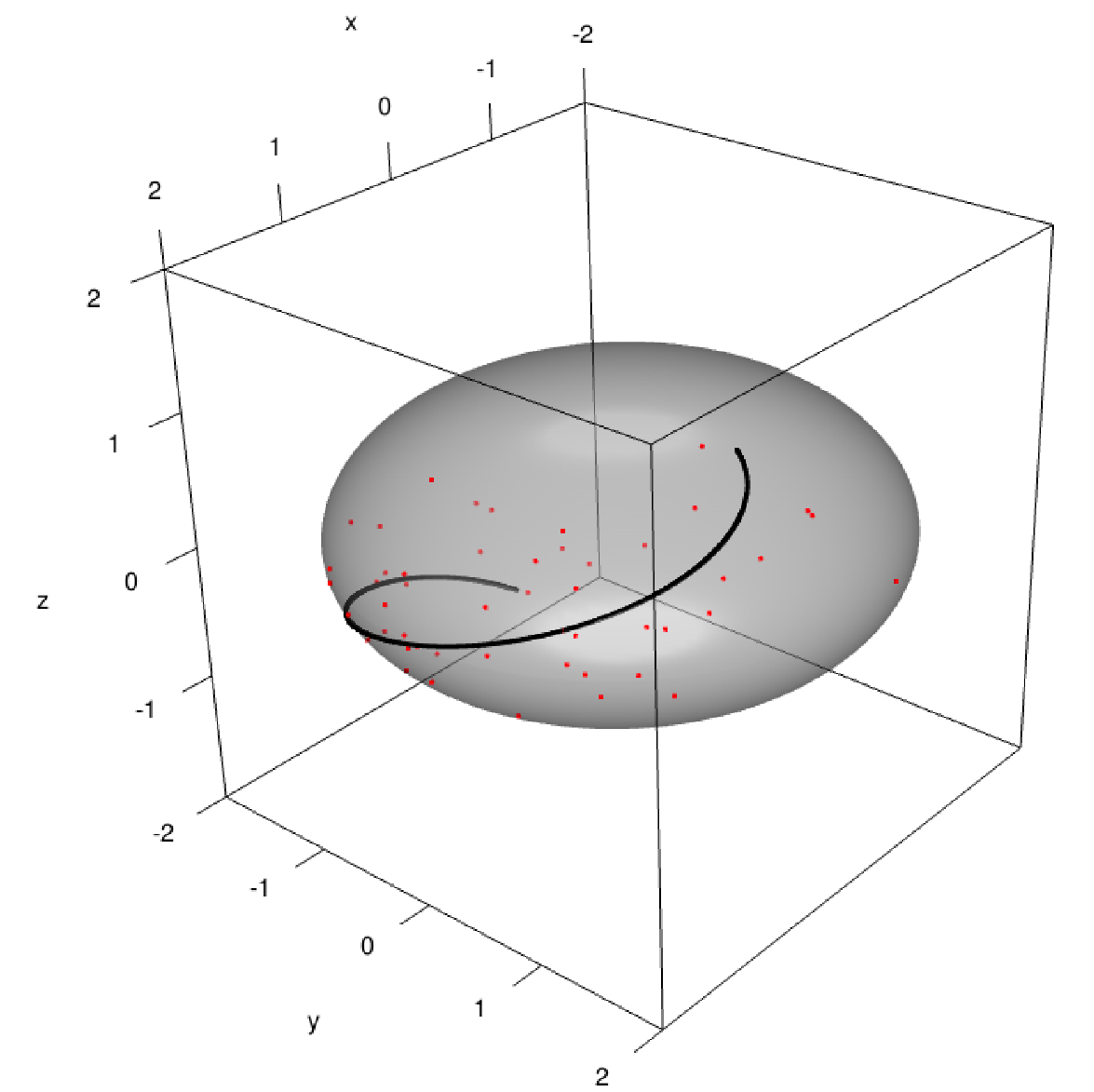}
    \subcaption{High Noise, $n = 50$}
  \end{minipage}
  \caption{\textit{Sample simulation data sets of size $n = 50$ under low (left) and high (right) noise settings.
The true regression curve is shown by the solid line.}}
\label{fig2}
\end{figure}

We generate a random sample \( (X_i, Y_i) \), \( i = 1, \ldots, n \), according to the following procedure.
First, we independently generate \( X_i \sim \mathcal{U}(-0.5, 0.5) \), for \( i = 1, \ldots, n \).
Next, we define a mapping \( s(x) \in S^2 \) by
\[
s(x) = \left( \sqrt{1 - 3x^2} \cos(\pi x),\ \sqrt{1 - 3x^2} \sin(\pi x),\ \sqrt{3} x \right).
\]
For each \( X_i \), we generate a zero-mean bivariate normal random vector $U_i$ with independent components
on the tangent space $T_{s(X_i)} S^2 = \left\{ v \in \mathbb{R}^3 \mid \langle v,\ s(X_i) \rangle = 0 \right\}$.
We then map this tangent vector to the sphere using the exponential map:
\[
Z_i = \mathrm{Exp}_{s(X_i)}(U_i) = \cos(\| U_i \|)\ s(X_i) + \sin(\| U_i \|)\ \frac{U_i}{\| U_i \|},
\]
where \( \| \cdot \| \) denotes the Euclidean norm and \( \mathrm{Exp} \) is the Riemannian exponential map (c.f., Example 5.4.1 in [\cite{MR2364186}]).
Finally, we obtain \( Y_i \in E \) by scaling the components of \( Z_i \) as follows: we multiply the \( x \)- and \( y \)-coordinates by 2 and leave the \( z \)-coordinate unchanged.
We conduct simulations with sample sizes \( n = 50,\ 100,\ 200 \), each repeated 100 times. We consider two noise levels for the variance of \( U_i \): 0.15 (low noise) and 0.30 (high noise). Fig.~\ref{fig2} illustrates simulated data sets with sample size $n = 50$.

We compare the LFR \eqref{2.5}, the LFR$^k$ \eqref{4.4}, and the Nadaraya-Watson (NW) smoother \eqref{3.3}.  
For estimation, a grid of bandwidths \( h \in [0.01, 0.15] \) is used for the smoothing, with $K$ being the 
Epanechnikov kernel: \( K(u) = \frac{3}{4}(1 - u^2)\mathbf{1}_{\{|u| \leq 1\}} \).  
For the LFR$^k$, the positive definite kernel \( k \) is chosen to be the Gaussian kernel  
$k(x, y) = \exp\left( - \frac{1}{\sigma^2} \|x - y\|_E^2 \right), \  \sigma = 1.5.$
For each sample size, we conduct 100 simulations and compute the mean integrated squared error (MISE) defined by  
\[
\mathrm{MISE}_r = \int_{-0.5}^{0.5} d^2\left( \hat{m}_{\oplus}^r(x), m_{\oplus}(x) \right) dx, \quad (r = 1, \ldots, 100),
\]
over a grid of bandwidths \( h \in [0.01, 0.15] \) with step size 0.01.  
Table~\ref{table2} reports the minimum average MISE across bandwidths and the corresponding optimal bandwidth.  
Table~\ref{table3} summarizes the average computation time per simulation for the case of \( n = 100 \).

\begin{table}
\centering
\caption{Minimum average MISE values (multiplied by 100 for clarity) for LFR, LFR$^k$, and NW estimators.  
The values in parentheses indicate the corresponding optimal bandwidth \( h \).}
\label{table2}
\begin{tabular}{lcccc}
\textbf{Noise}  \quad & \quad \ \ \textbf{n} \quad \ \ &
(a) \textbf{NW} \quad & \ \  \ \ \ \ (b) \textbf{LFR$^{k}$} \quad \  & (c) \textbf{LFR} \  \\
\hline \hline
Low   & 50  & 1.79(0.07) & 1.52(0.09) & 1.25(0.10) \\
    & 100 & 1.32(0.04) & 1.01(0.06) &  0.97(0.07)\\
    & 200 & 1.07(0.02) & 0.84(0.04)  & 0.80(0.04) \\ \hline
High   & 50  & 7.93(0.10) & 5.91(0.13)  & 5.62(0.14)\\
    & 100 & 6.53(0.06) & 4.08(0.10) & 3.89(0.11)\\
    & 200 & 4.55(0.04) & 2.66(0.06)  & 2.56(0.07) \\
\hline
\end{tabular}
\end{table}

\begin{table}
        \centering
          \caption{Average computation time per simulation (hours) for each estimator when \( n = 100 \).}
          \label{table3}
        \begin{tabular}{|c|c|c|c|}
            \hline
                       \ \ \textbf{Noise} \ \ & \ (a) \textbf{NW} \ & \ (b) \textbf{LFR$^{k}$} \ & \ (c) \textbf{LFR} \ \\ \hline  \hline
                                    Low & 1.60 & 0.076 & 1.61 \\ \hline
                                    High & 1.61 & 0.076 & 1.63 \\ \hline
                                    
        \end{tabular}
\end{table}

\begin{figure}
    \centering
    \begin{subfigure}{0.327\textwidth} 
        \centering
        \includegraphics[width=\linewidth]{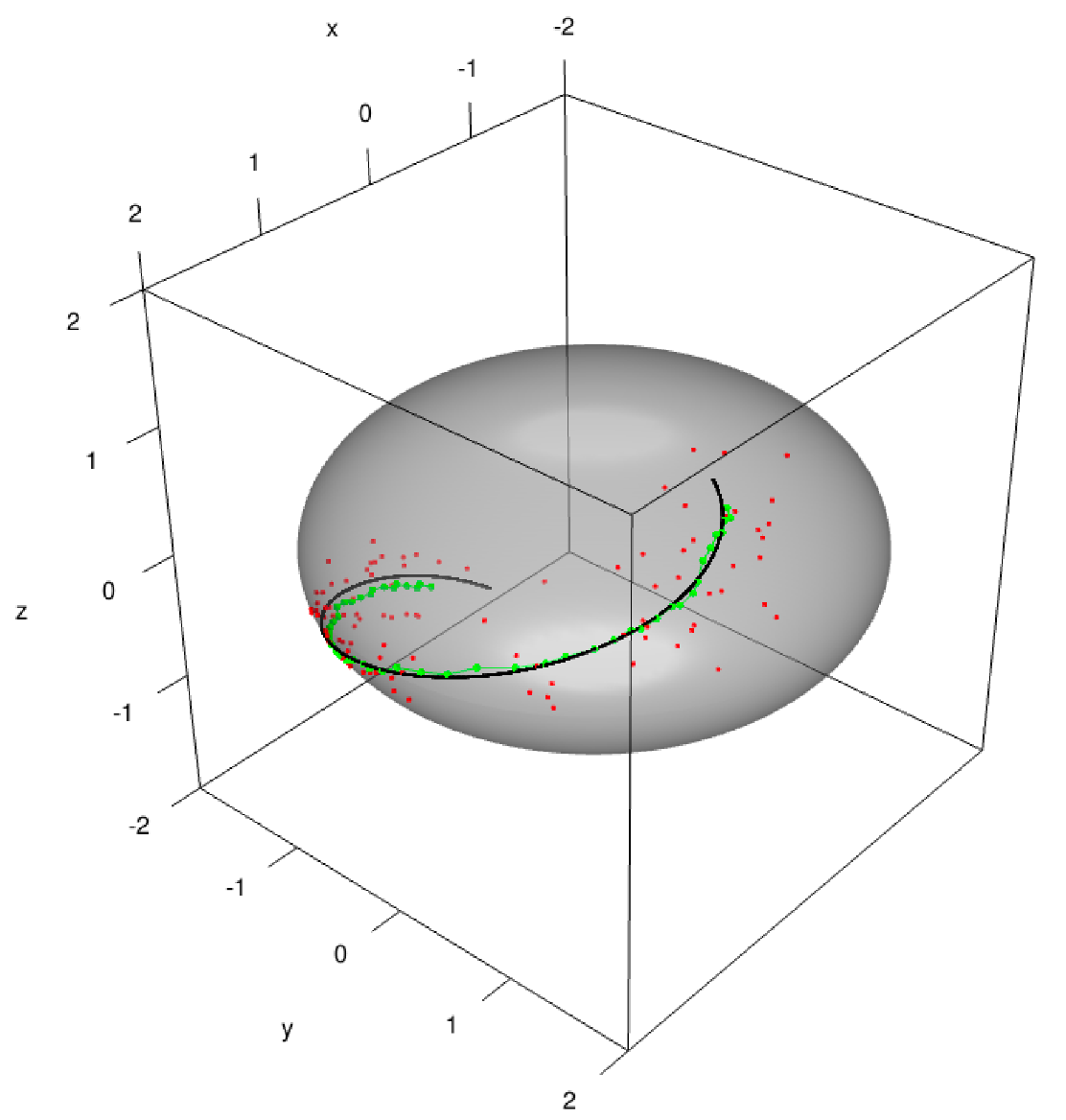} 
        \caption{NW} 
        \label{fig:a}
    \end{subfigure}
    \hfill     
    \begin{subfigure}{0.3281\textwidth}
        \centering
        \includegraphics[width=\linewidth]{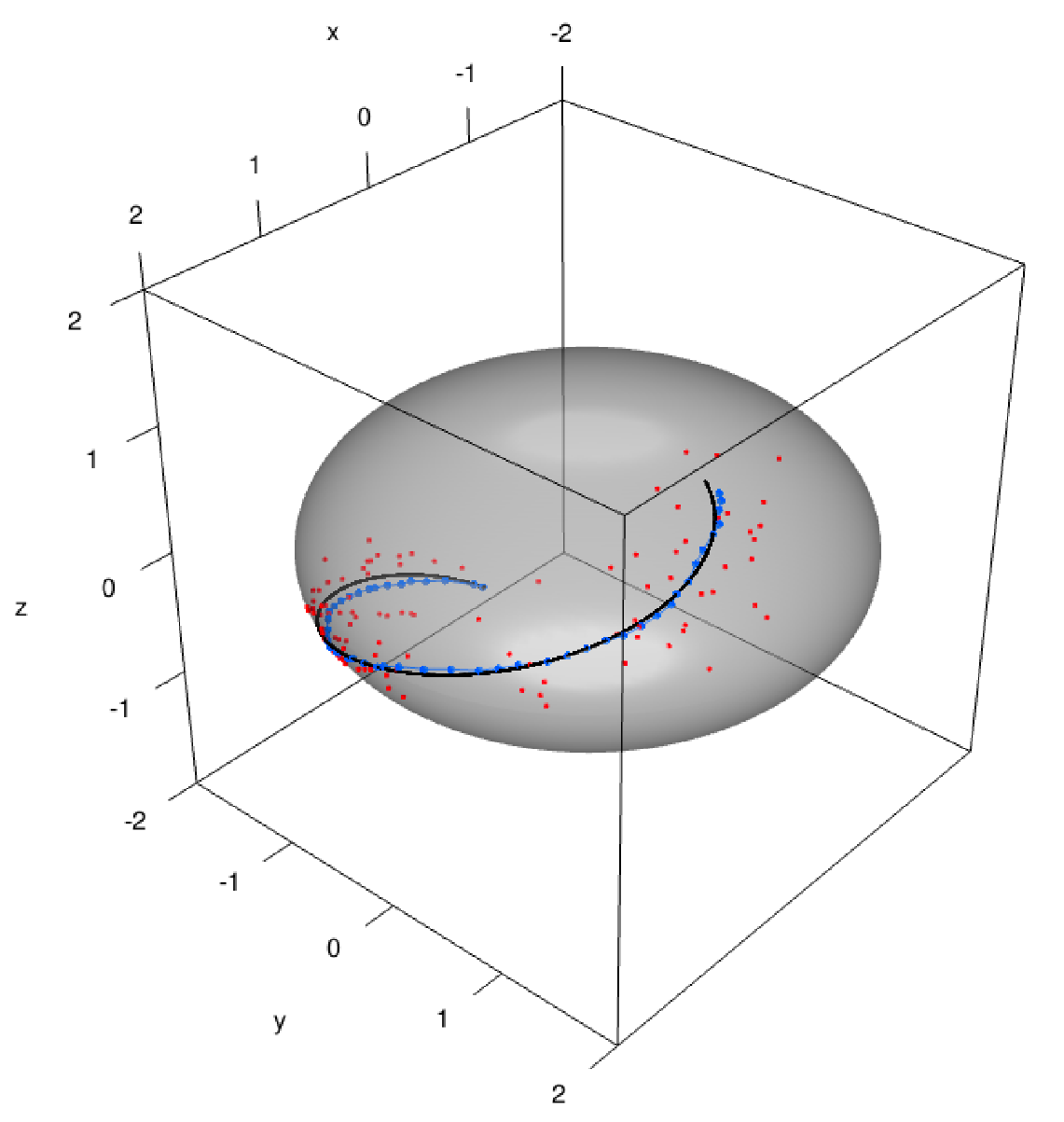}
        \caption{LFR$^k$}
        \label{fig:b}
    \end{subfigure}
    \hfill
    \begin{subfigure}{0.3262\textwidth}
        \centering
        \includegraphics[width=\linewidth]{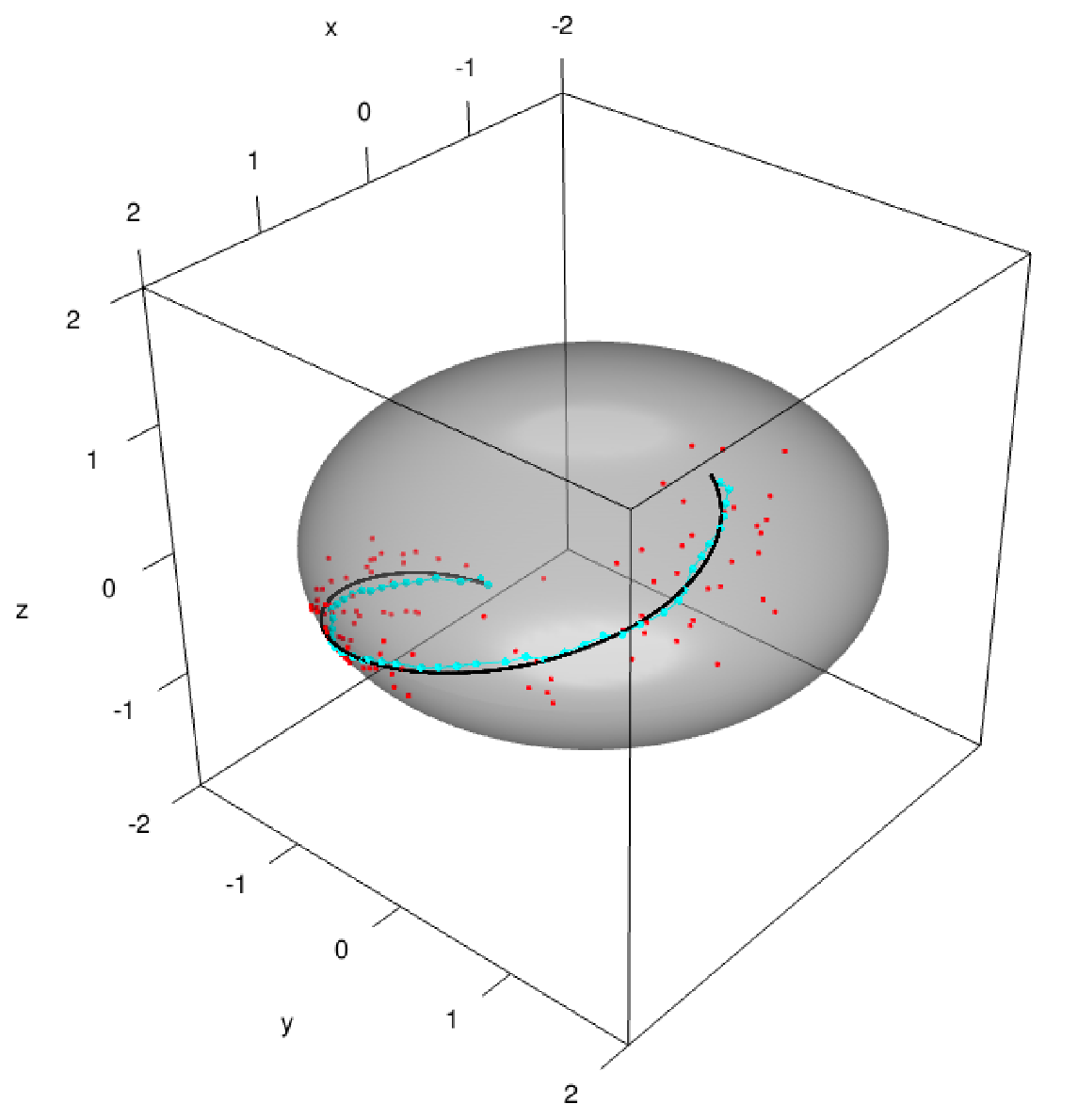}
        \caption{LFR}
        \label{fig:c}
    \end{subfigure}

    \caption{\textit{True regression curve $m_{\oplus}(x)$ (solid line) and the fitted values $\hat{m}_{\oplus}^{\mathrm{NW}}(x)$, 
    $\hat{l}_{\oplus}^k(x)$, and $\hat{l}_{\oplus}(x), \ x \in (-0.5, 0.5)$. (Low noise, $n=100$)}}    
     \label{fig3}
\end{figure}

Table~\ref{table2} shows that the MISE is ranked in the order of LFR, LFR$^k$, and NW in all settings.  
For LFR and LFR$^k$, it is observed that their values become closer as the sample size $n$ increases.  
On the other hand, Table~\ref{table3} indicates that the computation time of LFR$^k$ is approximately one twentieth of that for LFR and NW.  
This is because, when solving the minimization problems \eqref{2.5} and \eqref{3.3} for LFR and NW, numerical calculation of geodesic distances between two points on $E$ is required, whereas in the objective function \eqref{4.4} for LFR$^k$, geodesic distances are not involved, and only Euclidean distances inside the positive definite kernel need to be computed.  
Figure~\ref{fig3} presents the fitted values $\hat{l}_{\oplus}(x)$, $\hat{l}_{\oplus}^k(x)$, and $\hat{m}_{\oplus}^{\mathrm{NW}}(x)$ of LFR, LFR$^k$, and NW respectively on a grid of $x \in (-0.5, 0.5)$ with intervals of $0.02$.  
Compared to NW, both LFR and LFR$^k$ provide more accurate estimates near the boundaries of the true regression curve.

\begin{figure}
\centering
\includegraphics[width = 8cm]{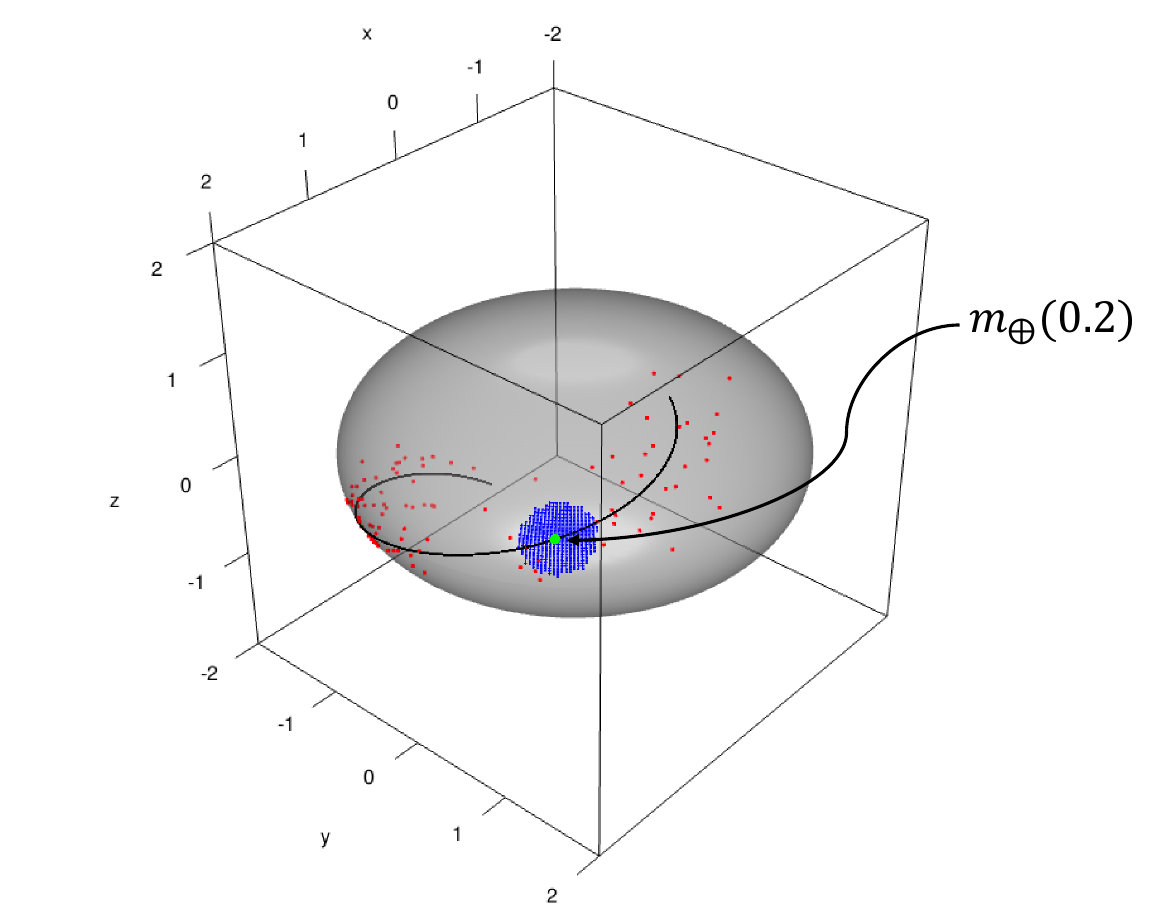} 
\caption{\textit{Approximate 95\% confidence region at the test point $x = 0.2$ using LFR$^k$. (Low noise, $n = 100$)}}
\label{fig4}
\end{figure}

Finally, based on the result of Corollary \ref{cr2}, we construct an approximate confidence region using the bootstrap method.  
An approximate confidence region can be obtained using the bias $\hat{B}_x$ and variance $\hat{V}_x$ estimated from bootstrap samples.  
Figure~\ref{fig4} illustrates the 95\% confidence region when the test point is $x = 0.2$.

\mbox{} \\

\begin{figure}
\centering
\includegraphics[width=5.9cm]{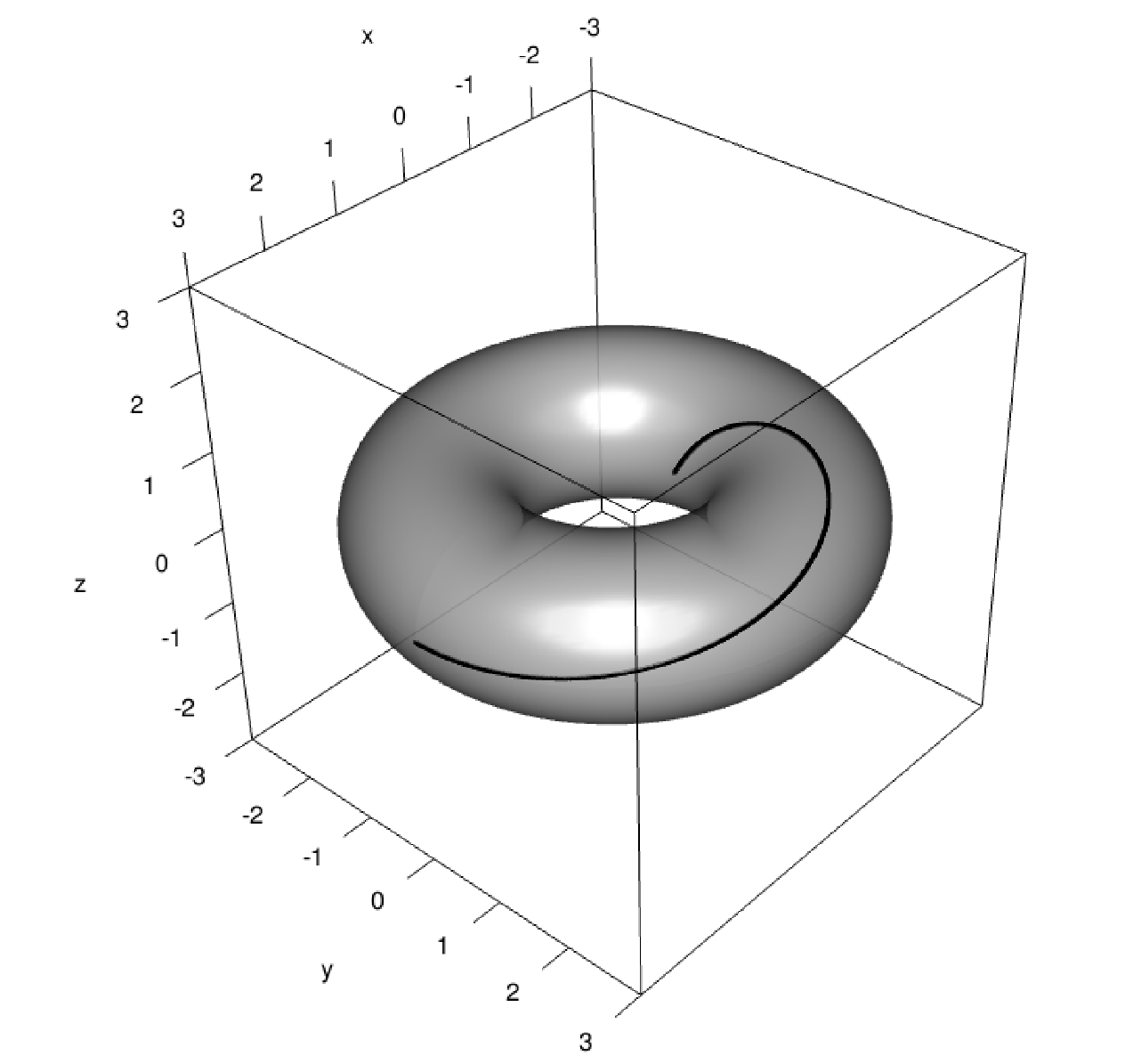}
\caption{\textit{True regression curve \( m_{\oplus}(x) \) on \( T \), \( x \in (0, 1) \).}}
\label{fig5}
\end{figure}

Next, we consider the torus \( T \) as the metric space \( \Omega \).  
A torus is a doughnut-shaped surface of revolution in three-dimensional space, generated by rotating a circle of radius \( r \) around an axis lying in the same plane as the circle and at a distance \( R > r \) from its center.  
The equation of a torus radially symmetric about the $z$-axis is given by 
\begin{align}
\Omega = T := \left\{ (x, y, z)^{\top} \in \mathbb{R}^3 \;\middle|\; \left( \sqrt{x^2 + y^2} - R \right)^2 + z^2 = r^2 \right\}. \label{5.3}
\end{align}
A point \( (x, y, z)^{\top} \) on \( T \) can be parameterized using two angular parameters \( \theta, \varphi \in [-\pi, \pi) \) as follows:
\begin{align*}
x = (R + r \cos \theta) \cos \varphi, \quad
y = (R + r \cos \theta) \sin \varphi, \quad
z = r \sin \theta.
\end{align*}
The distance function \( d \) on \( T \) is defined as the geodesic distance.
For two points \( P_1, P_2 \in T \), let \( \gamma: [0,1] \to T \) be any piecewise smooth curve connecting \( P_1 \) to \( P_2 \), and denote its length by \( L(\gamma) \).  
Then, the geodesic distance \( d \) on the torus is defined by
\begin{align*}
d(P_1, P_2) := \inf \left\{ L(\gamma) \;\middle|\; \gamma: [0,1] \to T,\ \gamma(0)=P_1,\ \gamma(1)=P_2 \right\},
\end{align*}
(c.f., Chapter 9 in  [\cite{spivak1970comprehensive}]).  
Except in the case of the standard torus defined by \( r = R \), this geodesic distance \( d \) does not admit a closed-form solution and must be computed numerically.  
Various methods have been proposed for numerically approximating geodesic distances on the torus, such as triangle mesh-based approximations [\cite{surazhsky2005fast}] and heat-based methods [\cite{crane2013geodesics}], though both are computationally intensive.
In the following, we conduct simulations on a torus \( T \) with parameters \( R = 2 \) and \( r = 1 \).
We consider the following true regression function on \( T \), as shown in Fig.~\ref{fig5}:
\begin{align}
m_{\oplus}(x) = \left( (2 + \cos(\pi x^2)) \cos(\pi x),\ (2 + \cos(\pi x^2)) \sin(\pi x),\ \sin(\pi x^2) \right),\ x \in (0, 1). \label{5.4}
\end{align}

\begin{figure}
   \centering
  \begin{minipage}[b]{0.44\linewidth}
    \centering
    \includegraphics[keepaspectratio, scale=0.26]{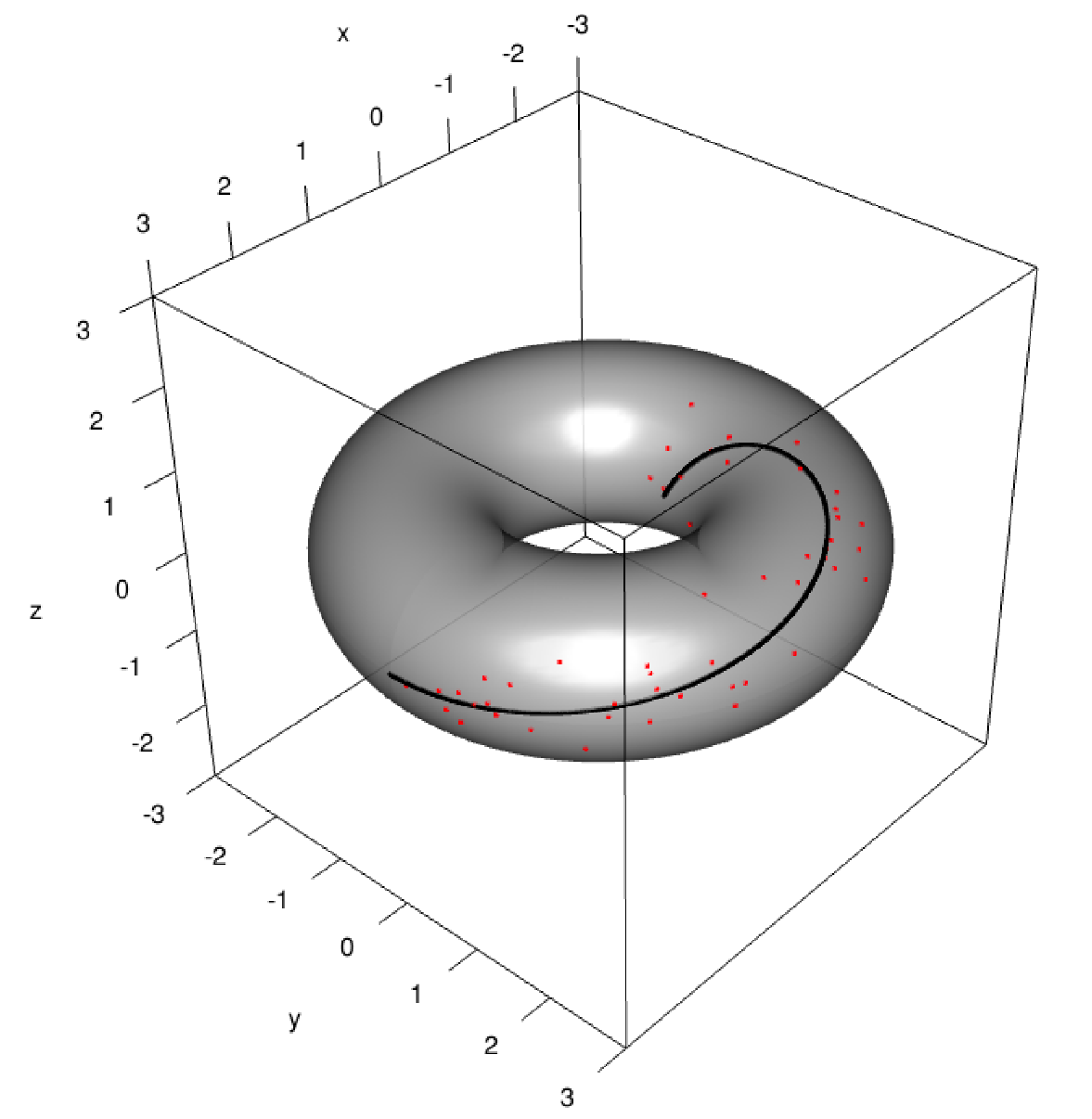}
    \subcaption{Low Noise, $n = 50$}
  \end{minipage}
  \begin{minipage}[b]{0.44\linewidth}
    \centering
    \includegraphics[keepaspectratio, scale=0.254]{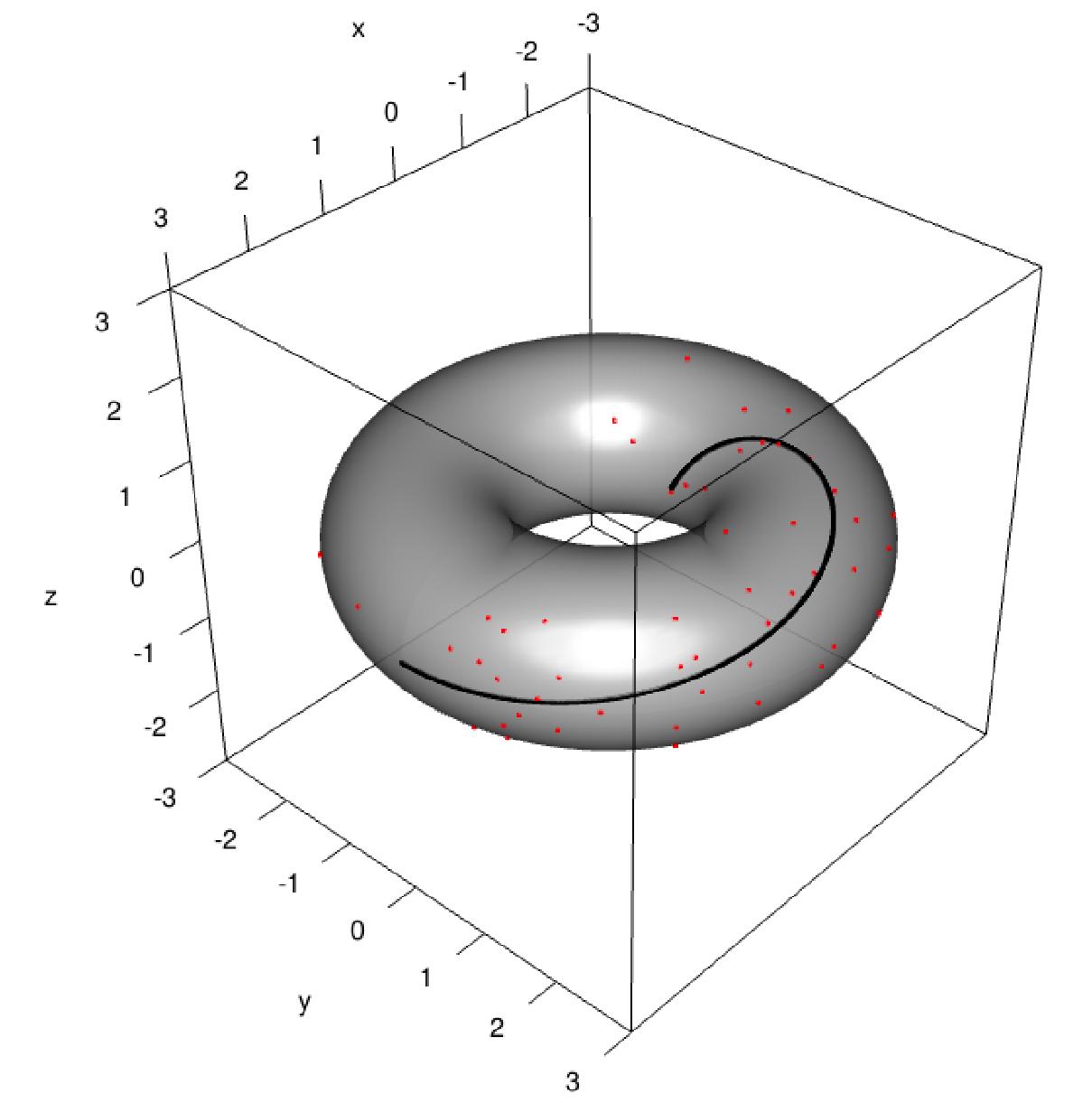}
    \subcaption{High Noise, $n = 50$}
  \end{minipage}
  \caption{\textit{Sample simulation data sets of size $n = 50$ under low (left) and high (right) noise settings.
The true regression curve is shown by the solid line.}}
\label{fig6}
\end{figure}

We generate a random sample \( (X_i, Y_i),\ i = 1, \ldots, n \) according to the following procedure.  
First, we independently generate \( X_i \sim \mathcal{U}(0, 1) \), for \( i = 1, \ldots, n \).
Then, for each \( X_i \), we generate two independent von Mises random variables \( (Y_{i1}, Y_{i2}) \) with mean directions \( \pi X_i \) and \( \pi X_i^2 \), respectively.  
The response \( Y_i \in T \) is defined by
\begin{align*}
Y_i = \left( (2 + \cos(Y_{i2})) \cos(Y_{i1}),\ (2 + \cos(Y_{i2})) \sin(Y_{i1}),\ \sin(Y_{i2}) \right).
\end{align*}
The von Mises distribution is a continuous probability distribution on the circle.  
Given a mean direction \( \mu \in [-\pi, \pi) \) and a concentration parameter \( \kappa \geq 0 \), its probability density function is given by (c.f., Section 3.5.4 in
[\cite{MR1828667}])
\begin{align*}
f_{\mathrm{vM}}(\theta;\mu,\kappa) = \frac{1}{2\pi I_0(\kappa)} \exp\left\{ \kappa \cos(\theta - \mu) \right\}, \quad \theta \in [-\pi, \pi),
\end{align*}
where \( I_0(\kappa) \) is the modified Bessel function of the first kind at order zero, defined as
\begin{align*}
I_0(\kappa) = \sum_{j=0}^\infty \frac{1}{(j!)^2} \left( \frac{\kappa}{2} \right)^{2j}.
\end{align*}
We conduct simulations for three sample sizes: \( n = 50, 100, 200 \).  
For each sample size, 100 independent simulations are performed.  
We consider two settings for the concentration parameter \( \kappa \):  
\( \kappa = 15 \) (low noise) and \( \kappa = 5 \) (high noise). 
Fig.~\ref{fig6} illustrates simulated data sets with sample size $n = 50$.

\begin{figure}
    \centering
    \begin{subfigure}{0.327\textwidth} 
        \centering
        \includegraphics[width=\linewidth]{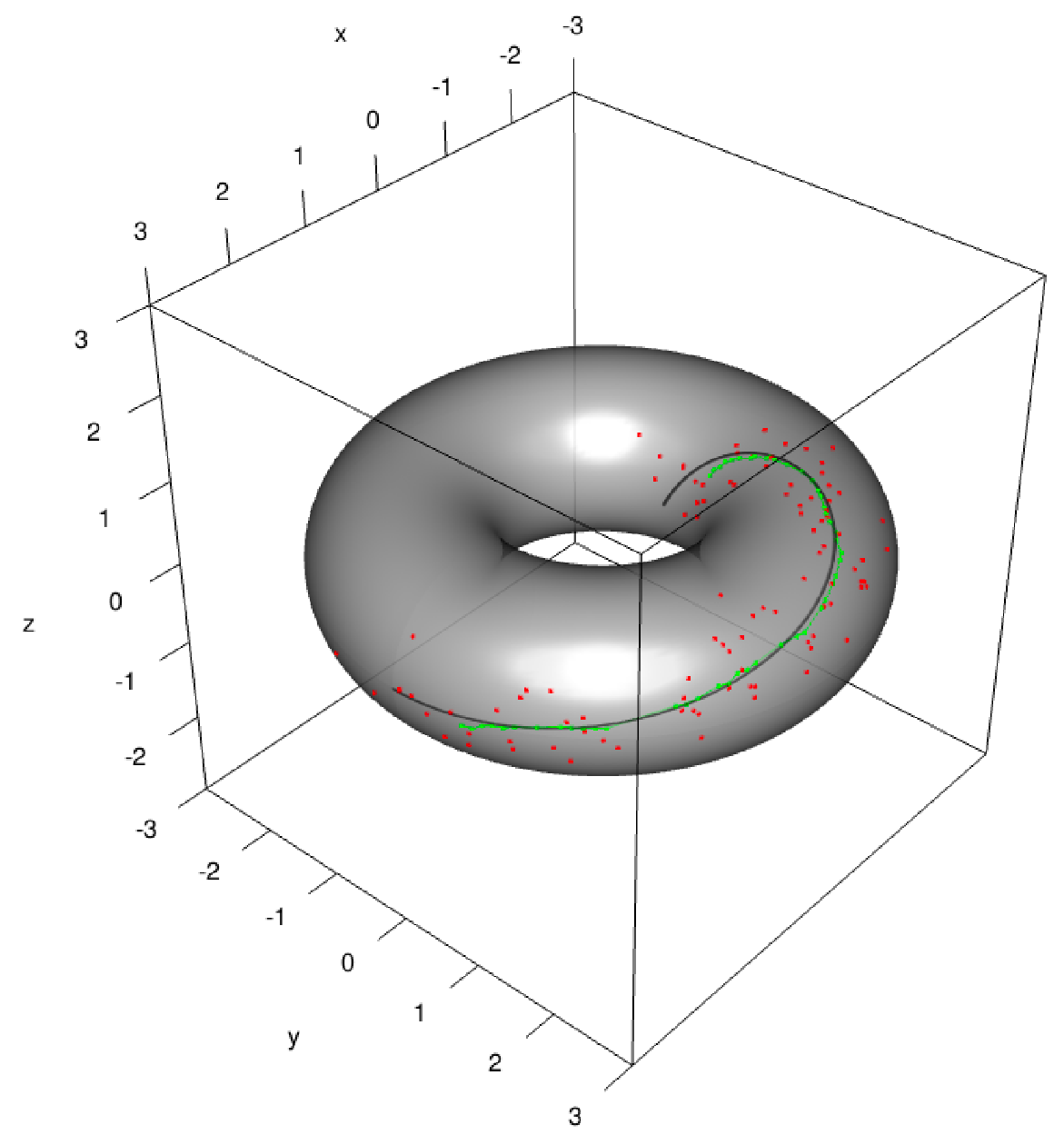} 
        \caption{NW}         \label{fig:a}
    \end{subfigure}
    \hfill 
    \begin{subfigure}{0.3285\textwidth}
        \centering
        \includegraphics[width=\linewidth]{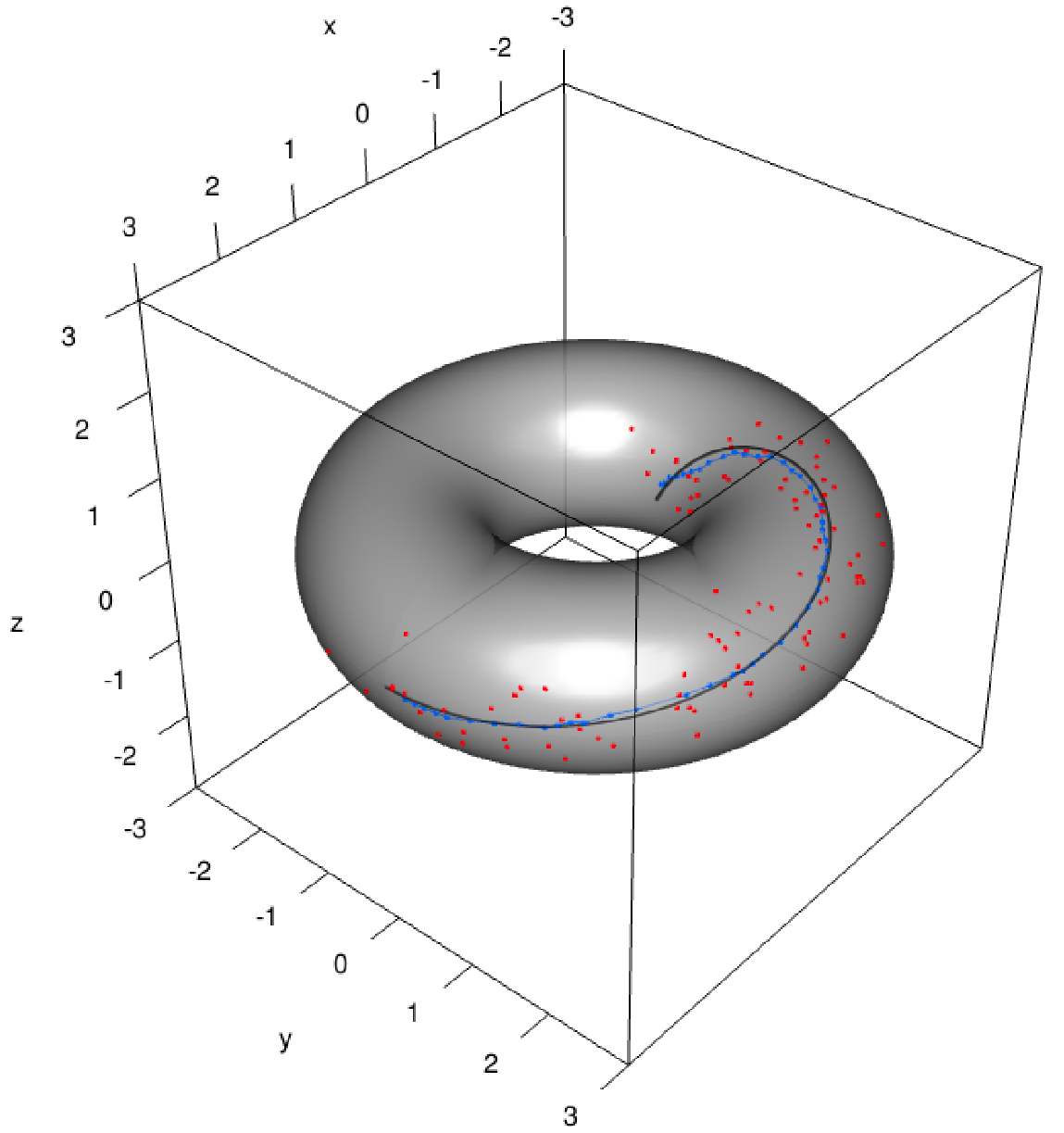}
        \caption{LFR$^k$}
        \label{fig:b}
    \end{subfigure}
    \hfill
    \begin{subfigure}{0.3262\textwidth}
        \centering
        \includegraphics[width=\linewidth]{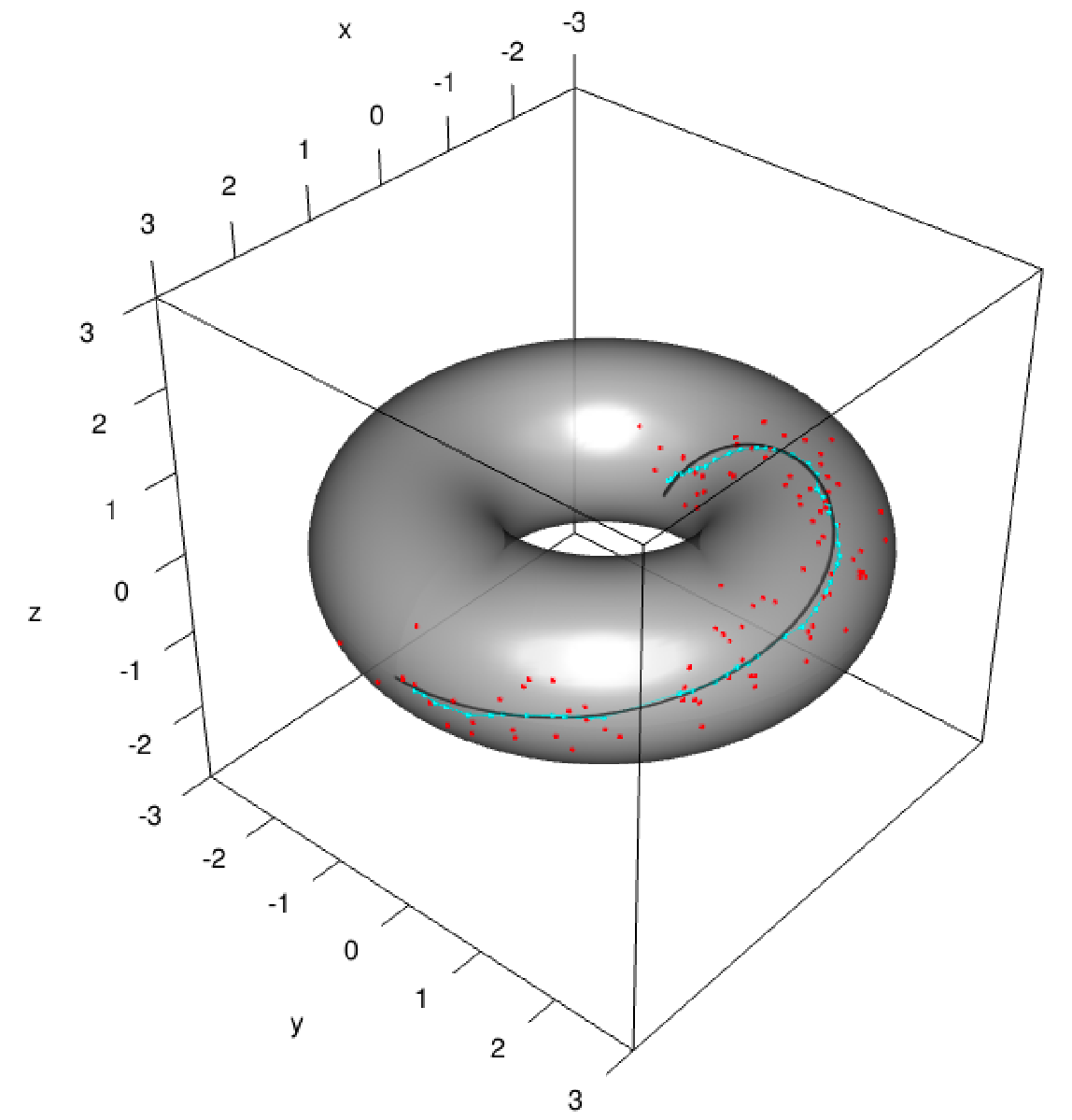}
        \caption{LFR}
        \label{fig:c}
    \end{subfigure}

    \caption{\textit{True regression curve $m_{\oplus}(x)$ (solid line) and the fitted values $\hat{m}_{\oplus}^{\mathrm{NW}}(x)$, 
    $\hat{l}_{\oplus}^k(x)$, and $\hat{l}_{\oplus}(x), \ x \in (0, 1)$. (Low noise, $n=100$)}}    
     \label{fig7}
\end{figure}

As in the previous example, we compare the LFR \eqref{2.5}, the LFR$^k$ \eqref{4.4}, and the Nadaraya-Watson (NW) smoother \eqref{3.3}.  
For estimation, a grid of bandwidths \( h \in [0.01, 0.15] \) is used for the smoothing, with $K$ being the 
Epanechnikov kernel: \( K(u) = \frac{3}{4}(1 - u^2)\mathbf{1}_{\{|u| \leq 1\}} \).  
For the LFR$^k$, the positive definite kernel \( k \) is chosen to be the Gaussian kernel  
$k(x, y) = \exp\left( - \frac{1}{\sigma^2} \|x - y\|_E^2 \right), \  \sigma = 1.5.$
The minimum values of the average MISE are reported in Table~\ref{table4}.
The average computation time per simulation for the case \( n = 100 \) is shown in Table~\ref{table5}.
Figure~\ref{fig7} presents the fitted values $\hat{l}_{\oplus}(x)$, $\hat{l}_{\oplus}^k(x)$, and $\hat{m}_{\oplus}^{\mathrm{NW}}(x)$ of LFR, LFR$^k$, and NW respectively on a grid of $x \in (0, 1)$ with intervals of $0.02$.  
Compared to NW, both LFR and LFR$^k$ demonstrate superior performance, especially near the boundaries.

\begin{table}
\centering
\caption{Minimum average MISE values (multiplied by 100 for clarity) for LFR, LFR$^k$, and NW estimators.  
The values in parentheses indicate the corresponding optimal bandwidth \( h \).}
\label{table4}
\begin{tabular}{lcccc}
\textbf{Noise}  \quad & \quad \ \ \textbf{n} \quad \ \ &
(a) \textbf{NW} \quad & \ \  \ \ \ \ (b) \textbf{LFR$^{k}$} \quad \  & (c) \textbf{LFR} \  \\
\hline \hline
Low   & 50  & 2.03(0.09) & 1.78(0.12) & 1.65(0.13) \\
      & 100 & 1.48(0.08) & 1.07(0.08) & 0.98(0.09) \\
      & 200 & 1.09(0.04) & 0.84(0.05) & 0.82(0.05) \\
\hline
High  & 50  & 8.26(0.11) & 5.72(0.13) & 5.30(0.13) \\
      & 100 & 6.29(0.09) & 4.29(0.11) & 4.07(0.12) \\
      & 200 & 4.55(0.06) & 2.76(0.08) & 2.70(0.09) \\
\hline
\end{tabular}
\end{table}

\begin{table}
        \centering
          \caption{Average computation time per simulation (hours) for each estimator when \( n = 100 \).}
          \label{table5}
        \begin{tabular}{|c|c|c|c|}
            \hline
                       \ \ \textbf{Noise} \ \ & \ (a) \textbf{NW} \ & \ (b) \textbf{LFR$^{k}$} \ & \ (c) \textbf{LFR} \ \\ \hline  \hline
                                    Low & 2.11 & 0.071 & 2.15 \\ \hline
                                    High & 2.13 & 0.072& 2.16 \\ \hline
                                    
        \end{tabular}
\end{table}

Through the above two simulation studies, it has been demonstrated that the proposed methods LFR and LFR$^k$ can accurately estimate regression curves even when the response variables lie on a manifold, where conventional regression methods for real-valued responses are not directly applicable. 
Moreover, by using LFR$^k$, it becomes unnecessary to compute geodesic distances on the manifold, resulting in a substantial reduction in computational cost.
\section{Conclusion}\label{Sec6}

In this paper, we focused on Local Fr\'echet Regression (LFR), a nonparametric regression method for response variables taking values in a metric space, and developed its asymptotic theory in the case where the response space is a Hilbert space. Specifically, by leveraging the fact that the LFR estimator admits an explicit representation in Hilbert spaces, we derived asymptotic normality under the framework of local polynomial regression.

Furthermore, to address situations in which direct application of LFR to general metric spaces is challenging, we proposed a novel estimation method (the LFR$^k$ estimator) that embeds the metric space data into a Reproducing Kernel Hilbert Space (RKHS) using a positive definite kernel, applies LFR in the RKHS, and then maps the result back to the original metric space. This approach allows for consistent estimation while avoiding computationally expensive optimization in the original metric space.

Finally, we conducted simulation studies using data with periodic structures lying on manifolds, which are difficult to handle using conventional statistical methods for real-valued responses. We demonstrated the effectiveness and computational efficiency of the proposed methods on random objects lying on a rotational ellipsoid and a torus, where the computation of geodesic distances is costly. In addition, we constructed an approximate confidence region for the regression function.

In summary, this work extends the theoretical foundation of LFR and introduces a practical estimation method that overcomes limitations in general metric spaces. It provides a valuable contribution to regression analysis in metric spaces and to the analysis of data lying on manifolds.

\bibliographystyle{imsart-nameyear} 
\bibliography{LFR_on_RKHS}
\begin{appendix}
\section*{PROOFS}\label{app} 
\subsection{Propositions \ref{pr1}, \ref{pr2}}\mbox{}\\

Since $\mathrm{E}[\| s(X, x, h)Y \|_\Omega] = \mathrm{E}[  |s(X, x, h)| \| Y \|_\Omega]
\leq  \mathrm{E}[ s(X, x, h)^2 ]^{1/2} \mathrm{E}[\| Y \|_{\Omega}^2 ]^{1/2} < \infty$,  
the map $\omega \mapsto \mathrm{E}[\langle s(X, x, h)Y, \omega \rangle _{\Omega}]$ is a bounded linear functional.  
Therefore, by the Riesz Representation Theorem (c.f., Chapter II, 16 in [\cite{MR1255973}]),  
there exists a unique $\tilde{\omega} \in \Omega$ such that for any $\omega \in \Omega$,  
\[
\mathrm{E}[\langle s(X, x, h)Y, \omega \rangle _{\Omega}] = \langle \tilde{\omega}, \omega \rangle _{\Omega}, 
\]
holds. Thus, noting that $\mathrm{E}[s(X, x, h)] = 1$, we obtain
\begin{align*}
\tilde{L}_n(\omega) &= \mathrm{E}[s(X, x, h) \| Y - \omega \|_\Omega^2]  \\
& = \mathrm{E}\left[ s(X, x, h) \left\{ \|Y - \tilde{\omega} \|_\Omega^2  + 2 \langle Y - \tilde{\omega} , \tilde{\omega}  - \omega  \rangle _\Omega       
+     \| \tilde{\omega}  - \omega  \| _\Omega^2    \right\} \right] \\
& = \mathrm{E}\left[ s(X, x, h)\|Y - \tilde{\omega}  \| _\Omega^2  \right] + 2 \mathrm{E}\left[ \langle s(X, x, h)Y , \tilde{\omega}  - \omega  \rangle _\Omega \right] \\
& \qquad \qquad \qquad \qquad - 2 \mathrm{E}\left[ s(X, x, h) \langle \tilde{\omega} , \tilde{\omega}  - \omega  \rangle _\Omega \right] 
+ \mathrm{E}\left[ s(X, x, h) \| \tilde{\omega}  - \omega  \|_\Omega^2 \right] \\
& = \mathrm{E}\left[ s(X, x, h)\|Y - \tilde{\omega}  \|_\Omega^2  \right] +  \|\tilde{\omega}  - \omega  \|_\Omega^2  .
\end{align*}
Since the first term on the right-hand side does not contain $\omega$, we have  
\[
\tilde{l}_{\oplus}(x) = \arg\min_{\omega \in \Omega} \tilde{L}_n(\omega) = \tilde{\omega}.
\]

Similarly, setting $\omega' = n^{-1} \sum_{i = 1}^{n} s_{in}(x, h) Y_i$ and noting that $n^{-1} \sum_{i = 1}^{n} s_{in}(x, h) = 1$,  
\begin{align*}
\hat{L}_n(\omega) &= n^{-1} \sum_{i = 1}^{n} s_{in}(x, h) \| Y_i - \omega \|_\Omega^2  \\
&= n^{-1} \sum_{i = 1}^{n} s_{in}(x, h) \left\{ \|Y_i - \omega' \|_\Omega^2  + 2 \langle Y_i - \omega' , \omega' - \omega  \rangle _\Omega        
+     \| \omega' - \omega  \| _\Omega^2    \right\} \\
& = n^{-1} \sum_{i = 1}^{n} s_{in}(x, h) \|Y_i - \omega' \| _\Omega^2 +  \| \omega' - \omega  \| _\Omega^2  .
\end{align*}
Since the first term on the right-hand side does not contain $\omega$, we have  
\[
\hat{l}_{\oplus}(x) = \arg\min_{\omega \in \Omega} \hat{L}_n(\omega) = \omega'.
\]
\qed
\subsection{Theorem \ref{th1}}\mbox{}\\

First, for the bivariate random sample $(X_1, Y_1), \ldots, (X_n, Y_n)$ from model \eqref{2.2}, the local linear estimator 
[\cite{MR1193323}, \cite{MR1212173}] $\hat{m}_{\mathrm{LL}}(x) = \hat{\beta}_0$ of $m(x)$ can be expressed using $s_{in}(x, h)$ as follows:
\begin{align*}
\hat{m}_{\mathrm{LL}}(x) = \sum_{i = 1}^{n} \frac{1}{n} \cdot \frac{\hat{\mu}_2 - \hat{\mu}_1 (X_i - x) }
{\hat{\mu}_2 \hat{\mu}_0 - \hat{\mu}_1 ^2} K_h(X_i - x) \cdot Y_i = n^{-1} \sum_{i = 1}^{n} s_{in}(x, h) Y_i.
\end{align*}
The estimator $\hat{m}_{\mathrm{LL}}(x)$ is asymptotically normal under (N1)--(N3) and the following (M1)--(M3) 
(c.f., Theorem 5.2 of [\cite{MR1383587}]).
\begin{description}
\item[(M1)]
$m(\cdot)$ is twice continuously differentiable on $N_x$. 

\item[(M2)]
$\sigma^2(\cdot)$ is continuous and positive on $N_x$.

\item[(M3)]
$\mathrm{E}[ Y^4 \mid X = \cdot ]$ is bounded on $N_x$.
\end{description}
\begin{align*}
\sqrt{nh} \left\{ \hat{m}_{\mathrm{LL}}(x) - m(x) - \left( \frac{ m''(x) }{2} \int u^2 K(u) du \right) h^2 \right\} 
 \leadsto \mathcal{N} \left( 0, \frac{\int K^2(u) du}{f(x)} \sigma^2 (x) \right).
\end{align*}

Next, let $\Omega$ be a real Hilbert space, and suppose a random sample $(X_1, Y_1), \ldots, (X_n, Y_n)$ is obtained from model \eqref{3.2}:
\begin{align*}
Y_i = m_{\oplus}(X_i) + \sigma(X_i) \varepsilon_i, \quad i = 1, \ldots, n.
\end{align*}
Taking the inner product of both sides with some $\omega \in \Omega$ yields:
\begin{align*}
\langle Y_i, \omega \rangle_{\Omega} = \langle m_{\oplus}(X_i), \omega \rangle_{\Omega} + \sigma(X_i) \langle \varepsilon_i, \omega \rangle_{\Omega}, \quad i = 1, \ldots, n.
\end{align*}
Define $Z_i = \langle Y_i, \omega \rangle_{\Omega},\  m_\omega(X_i) = \langle m_{\oplus}(X_i), \omega \rangle_{\Omega}$, and 
$\sigma(X_i) \sqrt{\mathrm{E}[\langle \varepsilon_i, \omega \rangle_{\Omega}^2]} \, \delta_i =  \sigma(X_i) \langle \varepsilon_i, \omega \rangle_{\Omega}$,
where $\delta_i = \langle \varepsilon_i, \omega \rangle_{\Omega} / \sqrt{\mathrm{E}[\langle \varepsilon_i, \omega \rangle_{\Omega}^2]}$, so that
$\mathrm{E}[\delta_i] = 0$, $\mathrm{Var}(\delta_i) = 1.$
Then we obtain the nonparametric regression model for $(X_1, Z_1), \ldots, (X_n, Z_n)$:
\begin{align*} 
Z_i = m_\omega(X_i) + \sigma(X_i) \sqrt{\mathrm{E}[\langle \varepsilon_i, \omega \rangle_{\Omega}^2]} \, \delta_i.  \tag{$\ast$} \label{A2}
\end{align*}
Using Proposition \ref{pr2}, the local linear estimator of $m_\omega(x)$ can be written using \( \hat{l}_\oplus(x) \) in \eqref{2.5} as:
\begin{align*}
\hat{m}_\omega^{\mathrm{LL}}(x) = n^{-1} \sum_{i = 1}^{n} s_{in}(x, h) Z_i
& = \left\langle n^{-1} \sum_{i = 1}^{n} s_{in}(x, h) Y_i, \ \omega \right\rangle_{\Omega} \\
& =  \left\langle \hat{l}_\oplus(x), \ \omega \right\rangle_{\Omega}
\end{align*}
The assumptions corresponding to (M1)--(M3) for model \hyperref[A2]{$(\ast)$} are denoted (O1)--(O3).
Applying the asymptotic normality of the local linear estimator to $\hat{m}_\omega^{\mathrm{LL}}(x) = \left\langle \hat{l}_\oplus(x), \ \omega \right\rangle_{\Omega}$,
we obtain that under (N1)--(N3) and (O1)--(O3), for any $\omega \in \Omega$,
\begin{align*}
&\sqrt{nh} \left\{ \hat{m}_\omega^{\mathrm{LL}}(x) - m_\omega(x) - \left( \frac{ m_\omega''(x) }{2} \int u^2 K(u) du \right) h^2 \right\} \\
{}= & \sqrt{nh} \left\{ \left\langle \hat{l}_\oplus(x) - m_\oplus(x), \ \omega \right\rangle_{\Omega} -  
 \left( \frac{\int u^2 K(u) du}{2}  \frac{\partial^2}{\partial x^2} \left\langle m_\oplus(x), \omega \right\rangle_{\Omega} \right) h^2 \right\} \\
{} \leadsto & \mathcal{N} \left( 0, \frac{\int K^2(u) du}{f(x)} \sigma^2(x) \mathrm{E}[\langle \varepsilon, \omega \rangle_{\Omega}^2] \right).
\end{align*}
\qed
\subsection{Corollary \ref{cr1}} \mbox{}\\

Let $Z \sim \mathcal{N}(0, 1)$. From the result of Theorem \ref{th1}, we have
\begin{align*}
\sqrt{nh} \, V_x^{-1/2} \left\{ \langle \hat{l}_{\oplus}(x), m_{\oplus}(x) \rangle_{\Omega} - \|m_{\oplus}(x)\|_{\Omega}^2 - B_x \right\} = Z + o_P (1).
\end{align*}
Suppose that consistent estimators $\hat{B}_x$ and $\hat{V}_x$ satisfy
\[
\hat{B}_x = B_x + o_ P \left((nh)^{-1/2}\right), \quad 
\hat{V}_x = V_x + o_P (1),
\]
then we obtain
\begin{align*}
\sqrt{nh} \, \hat{V}_x^{-1/2} \left\{ \langle \hat{l}_{\oplus}(x), m_{\oplus}(x) \rangle_{\Omega} - \|m_{\oplus}(x)\|_{\Omega}^2 - \hat{B}_x \right\} = Z + o_P (1).
\end{align*}
Therefore,
\begin{align*}
\langle \hat{l}_{\oplus}(x), m_{\oplus}(x) \rangle_{\Omega} 
= \|m_{\oplus}(x)\|_{\Omega}^2 + \hat{B}_x + \sqrt{\frac{\hat{V}_x}{nh}} \, Z + o_P \left((nh)^{-1/2}\right).
\end{align*}
Using this, we can derive the squared distance between the estimator and the true regression function as follows:
\begin{align*}
\|\hat{l}_{\oplus}(x) - m_{\oplus}(x)\|_{\Omega}^2 
&= \|\hat{l}_{\oplus}(x)\|_{\Omega}^2 + \|m_{\oplus}(x)\|_{\Omega}^2 - 2 \langle \hat{l}_{\oplus}(x), m_{\oplus}(x) \rangle_{\Omega} \\
&= \|\hat{l}_{\oplus}(x)\|_{\Omega}^2 - \|m_{\oplus}(x)\|_{\Omega}^2 - 2 \hat{B}_x - 2 \sqrt{\frac{\hat{V}_x}{nh}} \, Z + o_P \left((nh)^{-1/2}\right).
\end{align*}
According to Theorems 3 and 4 in [\cite{MR3909947}], under conditions (P1), (K0), and (L0)--(L3), we have
\begin{align*}
\left| \|\hat{l}_{\oplus}(x)\|_{\Omega}^2 - \|m_{\oplus}(x)\|_{\Omega}^2 \right|
&= \left| \langle \hat{l}_{\oplus}(x) + m_{\oplus}(x), \hat{l}_{\oplus}(x) - m_{\oplus}(x) \rangle_{\Omega} \right| \\
&\leq \|\hat{l}_{\oplus}(x) + m_{\oplus}(x)\|_{\Omega} \cdot \|\hat{l}_{\oplus}(x) - m_{\oplus}(x)\|_{\Omega} \\
&= O_P(1) \cdot O_P\left((nh)^{-1/(2(\beta_2 - 1))} + h^{2/(\beta_1 - 1)} \right) \\
&= o_P\left((nh)^{-1/2}\right).
\end{align*}
Hence, we obtain
\begin{align*}
\|\hat{l}_{\oplus}(x) - m_{\oplus}(x)\|_{\Omega}^2 
= -2 \hat{B}_x - 2 \sqrt{\frac{\hat{V}_x}{nh}} Z + o_P\left((nh)^{-1/2}\right),
\end{align*}
which leads to the approximation:
\begin{align*}
\sqrt{nh} \, \hat{V}_x^{-1/2} \left\{ \frac{1}{2} \|\hat{l}_{\oplus}(x) - m_{\oplus}(x)\|_{\Omega}^2 + \hat{B}_x \right\} = Z + o_P (1).
\end{align*}
Let $z_{\alpha/2}$ denote the upper $\alpha/2$ quantile of the standard normal distribution. Then, we obtain the asymptotic confidence statement:
\begin{align*}
P\left( 
- z_{\alpha/2} \leq 
\sqrt{nh} \, \hat{V}_x^{-1/2} \left\{ \frac{1}{2} \|\hat{l}_{\oplus}(x) - m_{\oplus}(x)\|_{\Omega}^2 + \hat{B}_x \right\}
\leq z_{\alpha/2} 
\right) = 1 - \alpha + o(1),
\end{align*}
which is equivalent to
\begin{align*}
P\left( 
-2\hat{B}_x - 2\sqrt{\frac{\hat{V}_x}{nh}} z_{\alpha/2} 
\leq \|\hat{l}_{\oplus}(x) - m_{\oplus}(x)\|_{\Omega}^2 
\leq -2\hat{B}_x + 2\sqrt{\frac{\hat{V}_x}{nh}} z_{\alpha/2}
\right) = 1 - \alpha + o(1).
\end{align*}
This establishes that
\[
P\left( m_{\oplus}(x) \in CI_{\alpha} \right) \to 1 - \alpha \quad \text{as } n \to \infty. 
\]
\qed
\subsection{Theorem \ref{th2}}\mbox{}\\

Given
$\hat{l}_{\oplus}^{H_k} (x) = n^{-1} \sum_{i = 1}^{n} s_{in}(x, h) \Phi(Y_i)$,\ 
$\tilde{l}_{\oplus}^{H_k} (x) = \argmin_{f \in H_k} \tilde{L}_n^{H_k}(f)$,
we define the corresponding empirical and population risk functions as follows:
\begin{align*}
 \hat{L}_n ^k(\omega) &= \|\hat{l}_{\oplus}^{H_k} (x) -  \Phi(\omega)\|^2 _{H_k} \\ 
& =  k(\omega, \omega)  - 2 {n}^{-1} \sum_{i = 1}^{n} s_{in}(x, h) k(Y_i, \omega) 
 + {n}^{-2} \sum_{i = 1}^{n}\sum_{j = 1}^{n} s_{in}(x, h) s_{jn}(x, h)k(Y_i, Y_j)  , \\
\tilde{L}_n ^k(\omega) &= \|\tilde{l}_{\oplus}^{H_k} (x) -  \Phi(\omega)\|^2 _{H_k} \\
 &= k(\omega, \omega) - 2\mathrm{E}[s(X, x, h) k(Y, \omega)] + \| \tilde{l}_{\oplus} ^{H_k}(x) \|_{H_k} ^2.
\end{align*}
Then, \eqref{4.4} can be written as 
$\hat{l}_{\oplus}^{k} (x) = \argmin_{\omega \in \Omega} \hat{L}_n ^k(\omega)$, 
and we similarly define \\$\tilde{l}_{\oplus}^{k}(x) = \argmin_{\omega \in \Omega} \tilde{L}_n^k(\omega)$. \\
By the triangle inequality for the distance function, it suffices to show
\begin{align*}
d(m_{\oplus}(x), \tilde{l}_{\oplus}^{k}(x)) = o(1) \quad &\cdots \text{(i)}, \\
d(\tilde{l}_{\oplus}^{k}(x), \hat{l}_{\oplus}^{k}(x)) = o_P(1) \quad &\cdots \text{(ii)}.
\end{align*}

First, regarding (i), by applying the same argument as in the proof of Theorem 3 in  SUPPLEMENTARY MATERIAL of [\cite{MR3909947}] 
(DOI:\href{https://doi.org/10.1214/17-AOS1624SUPP}{10.1214/17-AOS1624SUPP})
(hereafter referred to as [PM]),
and by replacing $\tilde{l}_{\oplus}(x)$ with $\tilde{l}_{\oplus}^{H_k}(x)$ and $m_{\oplus}(x)$ with $\Phi(m_{\oplus}(x))$,
under assumptions (K0), (L1), and (S4), we have
\[
\left\| \tilde{l}_{\oplus}^{H_k}(x) - \Phi(m_{\oplus}(x)) \right\|_{H_k}^2 = o(1) \quad \text{as } h = h_n \to 0.
\]
That is, since $\tilde{L}_n^k(m_{\oplus}(x)) \to 0$ as $h = h_n \to 0$, by (S5) it follows that
\[
d(m_{\oplus}(x), \tilde{l}_{\oplus}^{k}(x)) = o(1) \quad \text{as } h = h_n \to 0.
\]

Next, regarding (ii), we utilize Corollary 3.2.3 of [\cite{MR1385671}] (hereafter referred to as [VW]).  
Since assumption (S5) holds, it suffices to show
\[
\sup_{\omega \in \Omega} \left| \hat{L}_n^k(\omega) - \tilde{L}_n^k(\omega) \right| = o_P(1).
\]
To this end, we prove that
\[
\hat{L}_n^k - \tilde{L}_n^k \leadsto 0 \quad \text{in} \quad l^{\infty}(\Omega),
\]
and apply Theorem 1.3.6 in [VW] with
\[
g(L) = \sup_{\omega \in \Omega} |L(\omega)|, \quad L \in l^{\infty}(\Omega),
\]
where $l^{\infty}(\Omega)$ denotes the space of all uniformly bounded real-valued functions on $\Omega$.
By using Theorem 1.5.4 in [VW], this weak convergence is equivalent to showing the following two points:
\begin{itemize}
\item[(ii-\textcircled{\scriptsize 1})] $\hat{L}_n^k(\omega) - \tilde{L}_n^k(\omega) = o_P(1)$ for all $\omega \in \Omega$,
\item[(ii-\textcircled{\scriptsize 2})] $\hat{L}_n^k - \tilde{L}_n^k$ is asymptotically tight.
\end{itemize}
\quad
First, for (ii-\textcircled{\scriptsize 1}), the difference $\hat{L}_n^k(\omega) - \tilde{L}_n^k(\omega)$ can be written as:
\begin{align*}
&-2 \left\{ n^{-1} \sum_{i=1}^n s_{in}(x,h) k(Y_i,\omega) - \mathrm{E}[s(X,x,h)k(Y,\omega)] \right\} \\
&\quad \quad \quad \quad \quad+ n^{-2} \sum_{i=1}^n \sum_{j=1}^n s_{in}(x,h) s_{jn}(x,h) k(Y_i,Y_j) 
 - \| \tilde{l}_\oplus^{H_k}(x) \|_{H_k}^2.
\end{align*}
For the first term, by the result in the proof of Lemma 2 in [PM], it is $O_P((nh)^{-1/2})$. \\
For the second term, letting $s_i(x,h) = \sigma_0^{-2} K_h(X_i - x) \left\{ \mu_2 - \mu_1 (X_i - x) \right\}$,
we can write:
\begin{align*}
&n^{-2} \sum_{i=1}^n \sum_{j=1}^n \left\{ s_{in}(x,h) s_{jn}(x,h) - s_i(x,h) s_j(x,h) \right\} k(Y_i, Y_j) \tag{A-1} \label{A-1} \\
&\quad \quad \quad + n^{-2} \sum_{i=1}^n \sum_{j=1}^n s_i(x,h) s_j(x,h) k(Y_i, Y_j) - \| \tilde{l}_\oplus^{H_k}(x) \|_{H_k}^2 
\end{align*}
Now, for the first term in \eqref{A-1}, we compute the difference
$s_{in}(x,h) s_{jn}(x,h) - s_i(x,h) s_j(x,h)$ as:
\begin{align*}
&\left\{ \left (\frac{\hat{\mu}_2 }{\hat{\sigma}_0 ^2} \right)^2 -\left(\frac{\mu_2 }{\sigma_0 ^2} \right)^2 \right\} K_h(X_i - x) K_h(X_j - x) \\
& \quad\quad \quad - \left\{ \left(\frac{\hat{\mu}_1 }{\hat{\sigma}_0 ^2} \right)  \cdot  \left(\frac{\hat{\mu}_2 }{\hat{\sigma}_0 ^2} \right)
 - \left( \frac{\mu_1 }{\sigma_0 ^2} \right) \cdot \left( \frac{\mu_2 }{\sigma_0 ^2} \right) \right\} K_h(X_i - x) K_h(X_j - x) (X_i - x)\\
& \quad\quad \quad - \left\{ \left(\frac{\hat{\mu}_1 }{\hat{\sigma}_0 ^2} \right)  \cdot  \left(\frac{\hat{\mu}_2 }{\hat{\sigma}_0 ^2} \right)
 - \left( \frac{\mu_1 }{\sigma_0 ^2} \right) \cdot \left( \frac{\mu_2 }{\sigma_0 ^2} \right) \right\} K_h(X_i - x) K_h(X_j - x) (X_j - x)\\
&\quad \quad \quad \quad \quad \quad + \left\{ \left (\frac{\hat{\mu}_1 }{\hat{\sigma}_0 ^2} \right)^2 -\left(\frac{\mu_1 }{\sigma_0 ^2} \right)^2 \right\} 
 K_h(X_i - x) K_h(X_j - x) (X_i - x)(X_j - x).
\end{align*}
Lemma 1 in [PM], under assumptions (K0) and (L1), we have
\begin{align*}
\frac{\hat{\mu}_2 }{\hat{\sigma}_0 ^2}  - \frac{\mu_2}{\sigma_0 ^2} = O_{P}((nh)^{-1/2}), \quad
\frac{\hat{\mu}_1 }{\hat{\sigma}_0 ^2}  - \frac{\mu_1 }{\sigma_0 ^2} = O_{P}((nh^3)^{-1/2}).
\end{align*}
Therefore, it follows that
\begin{align*}
& \left( \frac{\hat{\mu}_2 }{\hat{\sigma}_0 ^2} \right)^2 - \left( \frac{\mu_2 }{\sigma_0 ^2} \right)^2 = O_{P}((nh)^{-1/2}), \quad
\left( \frac{\hat{\mu}_1 }{\hat{\sigma}_0 ^2} \right)^2 - \left( \frac{\mu_1 }{\sigma_0 ^2} \right)^2 = O_{P}((nh^3)^{-1/2}), \\
& \left( \frac{\hat{\mu}_1 }{\hat{\sigma}_0 ^2} \right) \cdot \left( \frac{\hat{\mu}_2 }{\hat{\sigma}_0 ^2} \right)
 - \left( \frac{\mu_1 }{\sigma_0 ^2} \right) \cdot \left( \frac{\mu_2 }{\sigma_0 ^2} \right) = O_{P}((nh^3)^{-1/2}).
\end{align*}
Moreover, under assumption (K0), we have
\begin{align*}
\mathrm{E}[K_h^i (X - x) (X - x)^j] 
= h^{j - i + 1} \int K^i(u) u^j f(x + hu) \, du 
= O(h^{j - i + 1}), \quad \text{for } i, j = 0,1,2,3,4 \ (i \neq 0).
\end{align*}
Therefore,
\begin{align*}
& K_h^2 (X_i - x) = O_{P}(h^{-3/2}), \quad K_h(X_i - x) K_h(X_j - x) = O_{P}(h^{-1}) \ (i \neq j), \\
& K_h^2 (X_i - x)(X_i - x) = O_{P}(h^{-1/2}), \quad K_h(X_i - x)K_h(X_j - x)(X_i - x) = O_{P}(1) \ (i \neq j), \\
& K_h^2 (X_i - x)(X_i - x)^2 = O_{P}(h^{1/2}), \quad K_h(X_i - x)K_h(X_j - x)(X_i - x)(X_j - x) = O_{P}(h) \ (i \neq j).
\end{align*}
Hence, we obtain
\begin{align*}
n^{-2} \sum_{i = 1}^{n} \sum_{j = 1}^{n} \left\{ s_{in}(x, h) s_{jn}(x, h) - s_i(x, h) s_j(x, h) \right\} k(Y_i, Y_j) = O_{P}((nh^4)^{-1/2}),
\end{align*}
and thus, under the condition $nh^4 \to \infty$, this term is $o_P(1)$. \\
Next, the second term of \eqref{A-1} can be written as:
\begin{align}
& \left\| n^{-1} \sum_{i = 1}^{n} s_i(x, h) \Phi(Y_i) - \tilde{l}_{\oplus}^{H_k}(x) \right\|_{H_k}^2  \tag{A-2} \label{A-2} \\
& \qquad \qquad \qquad \qquad 
+ 2 \left\langle \tilde{l}_{\oplus}^{H_k}(x), 
n^{-1} \sum_{i = 1}^{n} s_i(x, h) \Phi(Y_i) - \tilde{l}_{\oplus}^{H_k}(x) \right\rangle_{H_k}. \nonumber
\end{align}
Furthermore, define a functional $G: H_k \to \mathbb{R}$ for any $f \in H_k$ as
\begin{align*}
G(f) := \| f \|_{H_k}^2 + 2 \left\langle \tilde{l}_{\oplus}^{H_k}(x), f \right\rangle_{H_k}.
\end{align*}
Then $G$ is continuous, and expression \eqref{A-2} can be rewritten as
$G\left(n^{-1} \sum_{i = 1}^{n} s_i(x, h) \Phi(Y_i) - \tilde{l}_{\oplus}^{H_k}(x)\right)$.
From a result in the proof of Lemma 2 in [PM], for any $f \in H_k$, we have
\begin{align*}
\left\langle n^{-1} \sum_{i = 1}^{n} s_i(x, h) \Phi(Y_i) - \tilde{l}_{\oplus}^{H_k}(x), f \right\rangle_{H_k}
&= n^{-1} \sum_{i = 1}^{n} s_i(x, h) f(Y_i) - \mathbb{E}[s(X, x, h) f(Y)] \\
&= O_P((nh)^{-1/2}),
\end{align*}
and thus the left-hand side converges in probability to $0 = \langle 0, f \rangle$ as $nh \to \infty$.
Under assumption (S1), since $H_k$ is separable, there exists a countable orthonormal basis $\{e_i\}_{i=1}^\infty \subset H_k$. Then, for any $j \in \mathbb{N}$, the inner product between $e_j$ and $n^{-1} \sum_{i = 1}^{n} s_{in}(x, h) \Phi(Y_i) - \tilde{l}_{\oplus}^{H_k}(x)$ converges to zero in probability. Therefore, $n^{-1} \sum_{i = 1}^{n} s_{in}(x, h) \Phi(Y_i) - \tilde{l}_{\oplus}^{H_k}(x)$
is asymptotically finite-dimensional.
By Theorem 1.8.4 in [VW], we obtain
\begin{align*}
n^{-1} \sum_{i = 1}^{n} s_{in}(x, h) \Phi(Y_i) - \tilde{l}_{\oplus}^{H_k}(x) \leadsto 0.
\end{align*}
Applying Theorem 1.3.6 in [VW] to the functional $G$, the norm and inner product in \eqref{A-2} satisfy
$G\left(n^{-1} \sum_{i = 1}^{n} s_i(x, h) \Phi(Y_i) - \tilde{l}_{\oplus}^{H_k}(x)\right) = o_P(1)$.
Therefore, we conclude that for any $\omega \in \Omega$, $\hat{L}_n^k(\omega) - \tilde{L}_n^k(\omega) = o_P(1)$.

Next, regarding (ii-\textcircled{\scriptsize 2}), by applying Theorem 1.5.7 in [VW], it suffices to show that \\
$\hat{L}_n^k - \tilde{L}_n^k$ is asymptotically uniformly $d$-equicontinuous in probability, i.e., for any $\eta > 0$,
\begin{align*}
\limsup_{n} P \left( \sup_{d(\omega_1, \omega_2) < \delta}  
\left |  (\hat{L}_n^k  - \tilde{L}_n^k)(\omega_1) - (\hat{L}_n^k  - \tilde{L}_n^k)(\omega_2) \right | > \eta \right) \to 0 \quad \text{as } \delta \to 0. \tag{A-3} \label{A-3}
\end{align*}
Since $\mathrm{E}[|s(X, x, h)|] = O(1)$, $\mathrm{E}[s^2(X, x, h)] = O(h^{-1})$, and 
${n}^{-1} \sum_{i=1}^{n} | s_{in}(x, h) | = O_{P}(1)$, 
under assumptions (S2) and (S3), we have
\begin{align*}
\left | \hat{L}_n^k(\omega_1) - \hat{L}_n^k(\omega_2) \right | 
& \leq  \left | k(\omega_1, \omega_1) - k(\omega_2, \omega_2) \right | + 2 \left | \frac{1}{n} \sum_{i = 1}^{n} s_{in}(x, h) \left \{ k(Y_i, \omega_2) - k(Y_i, \omega_1) \right \} \right | \\
& \leq \left |  \| \Phi(\omega_2) \|^2_{\mathcal{H}_k} - \| \Phi(\omega_1) \|^2_{\mathcal{H}_k} \right | \\
& \quad \quad + \frac{1}{n} \sum_{i = 1}^{n} | s_{in}(x, h) | \times \\
& \quad \quad \quad \quad \left | \|\Phi(Y_i) - \Phi(\omega_1)\|^2_{\mathcal{H}_k} - \|\Phi(Y_i) - \Phi(\omega_2)\|^2_{\mathcal{H}_k} 
 + \| \Phi(\omega_2) \|^2_{\mathcal{H}_k} - \| \Phi(\omega_1) \|^2_{\mathcal{H}_k} \right | \\
& \leq 2MC\, d(\omega_1, \omega_2) + 
\frac{1}{n} \sum_{i = 1}^{n} | s_{in}(x, h) | \times 2C \left\{C\, \mathrm{diam}(\Omega)\, d(\omega_1, \omega_2) + M\, d(\omega_1, \omega_2) \right\} \\
& = O_{P}(d(\omega_1, \omega_2)).
\end{align*}
Here, we use the fact that, by assumption (S1),  
$d(\omega, \omega') \leq \mathrm{diam}(\Omega) < \infty$ holds for any $\omega, \omega' \in \Omega$.\\
Similarly, we have
\begin{align*}
\left| \tilde{L}_n^k(\omega_1) - \tilde{L}_n^k(\omega_2) \right| 
& \leq 2MC\, d(\omega_1, \omega_2) + \mathrm{E}[ | s(X, x, h) | ] \times 2C \left\{ C\,\mathrm{diam}(\Omega)\, d(\omega_1, \omega_2) + M\, d(\omega_1, \omega_2) \right\} \\
& = O( d(\omega_1, \omega_2) ).
\end{align*}
From the above, we conclude that for any $\eta > 0$, the condition \eqref{A-3} holds.
\subsection{Corollary \ref{cr2}} \mbox{}\\

Suppose that $(X, Y) \in \mathbb{R} \times \Omega$ is defined on the probability space $(T \times U, \mathcal{F}_T \times \mathcal{F}_U, P)$. 
Define the following set of functions in \( H_k \):
\begin{align*}
CI_{\alpha}^{H_k}
:= \left\{ g \in H_k \ \middle| \ 
-2 \hat{B}_x^{H_k} - 2 \sqrt{ \frac{\hat{V}_x^{H_k}}{nh} } z_{\alpha/2}
\le 
\| g - \hat{l}_{\oplus}^{H_k}(x) \|^2_{H_k}
\le 
-2 \hat{B}_x^{H_k} + 2 \sqrt{ \frac{\hat{V}_x^{H_k}}{nh} } z_{\alpha/2}
\right\}.
\end{align*}
Following the same reasoning as in Corollary \ref{cr1}, Theorem \ref{th3} yields
\[
P\left( \Phi(m_{\oplus}(x)) \in CI_{\alpha}^{H_k} \right) \to 1 - \alpha \quad \text{as } n \to \infty.
\]
For any \((t, u) \in T \times U\), the values \((X_i(t), Y_i(u))\) for \(i = 1, \ldots, n\) are determined, 
and based on these realizations, we obtain non-random (deterministic) confidence regions $CI_{\alpha}^{H_k}(t, u) \subset H_k$ and $CI_{\alpha}(t, u) \subset \Omega$.
From the definitions of $CI_{\alpha}^{H_k}$ and $CI_{\alpha}$, it follows that
\[
\Phi(m_{\oplus}(x)) \in CI_{\alpha}^{H_k}(t, u) \iff m_{\oplus}(x) \in CI_{\alpha}(t, u).
\]
Thus, we have
\[
P\left( m_{\oplus}(x) \in CI_{\alpha} \right) \to 1 - \alpha \quad \text{as } n \to \infty.
\]
Furthermore, Theorem \ref{th2} implies that
\[
P\left( \hat{l}_{\oplus}^{k}(x) \in CI_{\alpha} \right) \to 1 - \alpha \quad \text{as } n \to \infty. 
\]
\qed

\end{appendix}
\begin{supplement}\label{S1}
\stitle{S.1\ \ Definition \ref{de4}.}
\end{supplement}

Under the same model as in \eqref{3.2}, namely, $Y = \mathfrak{m}(X) + \sigma(X)\varepsilon$, we have
\begin{align*}
M_{\oplus}(\omega, x) &= \mathrm{E}\left[\| Y - \omega \|^2_{\Omega} \mid X = x\right] \\
&= \mathrm{E}\left[\| \mathfrak{m}(x) + \sigma(x)\varepsilon - \omega \|^2_{\Omega} \right] \\
&= \mathrm{E}\left[ \langle \sigma(x)\varepsilon + \mathfrak{m}(x) - \omega, \sigma(x)\varepsilon + \mathfrak{m}(x) - \omega \rangle_{\Omega} \right] \\
&= \sigma^2(x) \mathrm{E}\left[ \| \varepsilon \|^2_{\Omega} \right] 
+ 2\sigma(x)\mathrm{E}\left[ \langle \varepsilon, \mathfrak{m}(x) - \omega \rangle_{\Omega} \right] 
+ \| \mathfrak{m}(x) - \omega \|^2_{\Omega} \\
&= \sigma^2(x) + \| \mathfrak{m}(x) - \omega \|^2_{\Omega}.
\end{align*}
Therefore,
\[
m_{\oplus}(x) = \arg\min_{\omega \in \Omega} M_{\oplus}(\omega, x) = \mathfrak{m}(x).
\]

\begin{supplement}\label{S2}
\stitle{S.2\ \ Table \ref{table1}.} 
\end{supplement}

For the local linear estimator [\cite{MR1193323}, \cite{MR1212173}]
\( \hat{m}_{\mathrm{LL}}(x) \), the conditional bias and conditional variance are given as follows  
(c.f., Table 1 of [\cite{MR1209561}], Theorem 3.1 of [\cite{MR1383587}]):
Under assumptions (N1)--(N3), (M1), and (M2), 
\begin{align*}
\mathrm{E}[\hat{m}_{\mathrm{LL}}(x) \mid X_1, \ldots, X_n] - m(x) &= \frac{\int u^2 K(u) \, du}{2} \, m''(x) \, h^2 + o_P(h^2), \\
\mathrm{Var}(\hat{m}_{\mathrm{LL}}(x) \mid X_1, \ldots, X_n) &= \frac{\int K^2(u) \, du}{f(x) \, nh} \, \sigma^2(x) + o_P((nh)^{-1}).
\end{align*}
For the bivariate random sample \( (X_1, Y_1), \ldots, (X_n, Y_n) \) generated from the model \eqref{2.2},  
the Nadaraya-Watson estimator of \( m(x) \) ([\cite{MR166874}], [\cite{MR185765}]) is defined by
\[
\hat{m}_{\mathrm{NW}}(x) = \frac{\sum_{i=1}^n K_h(X_i - x) Y_i}{\sum_{j=1}^n K_h(X_j - x)}.
\]
For \( \hat{m}_{\mathrm{NW}}(x) \), the conditional bias and conditional variance are given as follows  
(c.f., Table 1 of [\cite{MR1209561}], Theorem 3.1 of [\cite{MR1383587}]):
Under assumptions (N1)--(N3), (M1) and (M2), further assume that \( f'(\cdot) \) is continuous on \( N_x \) and \( nh^3 \to \infty \), then
\begin{align*}
\mathrm{E}[\hat{m}_{\mathrm{NW}}(x) \mid X_1, \ldots, X_n] - m(x) 
&= \frac{\int u^2 K(u) \, du}{2} \left\{ m''(x) + 2 \frac{m'(x) f'(x)}{f(x)} \right\} h^2 + o_P(h^2), \\
\mathrm{Var}(\hat{m}_{\mathrm{NW}}(x) \mid X_1, \ldots, X_n) 
&= \frac{\int K^2(u) \, du}{f(x) nh} \, \sigma^2(x) + o_P((nh)^{-1}).
\end{align*}
The Nadaraya-Watson estimator for \( m_\omega(x) \) can be written using \( \hat{m}^{\mathrm{NW}}_{\oplus}(x) \) in \eqref{3.3} as:
\begin{align*}
\hat{m}^{\mathrm{NW}}_\omega(x) 
= \frac{\sum_{i = 1}^{n} K_h(X_i - x) Z_i}{\sum_{j = 1}^{n} K_h(X_j - x)} 
&= \left\langle \frac{\sum_{i = 1}^{n} K_h(X_i - x) Y_i}{\sum_{j = 1}^{n} K_h(X_j - x)}, \, \omega \right\rangle_{\Omega} \\
&= \left\langle \hat{m}^{\mathrm{NW}}_{\oplus}(x), \, \omega \right\rangle_{\Omega}.
\end{align*}
Therefore, applying the above results to  
\[
\hat{m}^{\mathrm{LL}}_\omega(x) = \left\langle \hat{l}_{\oplus}(x), \, \omega \right\rangle_{\Omega}
\quad \text{and} \quad
\hat{m}^{\mathrm{NW}}_\omega(x) = \left\langle \hat{m}^{\mathrm{NW}}_{\oplus}(x), \, \omega \right\rangle_{\Omega},
\]
we obtain the results summarized in Table \ref{table1}.

\end{document}